\documentclass[11pt]{amsart}
\usepackage{amsmath,amsthm, amscd, amssymb, amsfonts}
\usepackage[all]{xy}
\usepackage{hhline}
\usepackage[dvips, dvipsnames, usenames]{color}
\usepackage{eucal}

\newcommand{\ot}{\otimes}
\newcommand{\Ind}{\operatorname{Ind}}

\newcommand{\diag}{\operatorname{diag}}

\newcommand{\kc}{\mathbb F_q}
\newcommand{\Tb}{\mathbb T}

\newcommand{\gi}{g_{\infty}}

\newcommand{\oct}{\mathfrak O}
\newcommand{\q}{\mathfrak q}
\newcommand\toba{{\mathfrak B }}
\newcommand\trasp{\sigma}

\newcommand{\trid}{\triangleright}

\newcommand{\W}{{\mathcal W}}

\newcommand{\ku}{\mathbb C}

\newcommand{\Z}{{\mathbb Z}}
\newcommand{\N}{{\mathbb N}}

\newcommand{\G}{{\mathbb G}}

\newcommand{\D}{{\mathcal D}}

\newcommand{\Oc}{{\mathcal O}}
\newcommand{\oc}{{\mathcal O}}

\newcommand{\ydga}{{}^{\ku\,\Gamma}_{\ku\,\Gamma}\mathcal{YD}}
\newcommand{\ydg}{{}^{\ku\,G}_{\ku\,G}\mathcal{YD}}

\newcommand{\Aut}{\operatorname{Aut}}

\newcommand\sgn{\operatorname{sgn}}

\numberwithin{equation}{section}
\theoremstyle{plain}

\newtheorem{lema}{Lemma}[section]
\newtheorem{theorem}[lema]{Theorem}
\newtheorem{cor}[lema]{Corollary}

\newtheorem{prop}[lema]{Proposition}
\newtheorem{example}[lema]{Example}

\newtheorem{exa}[lema]{Example}

\theoremstyle{definition}
\newtheorem{definition}[lema]{Definition}

\theoremstyle{remark}
\newtheorem{obs}[lema]{Remark}
\newtheorem{remark}[lema]{Remark}
\newtheorem{rmk}[lema]{Remarks}

\newcommand\id{\operatorname{id}}

\newcommand\Id{\operatorname{id}}

\newcommand\s{\mathbb S}
\newcommand\st{\mathbb S_3}
\newcommand\sk{\mathbb S_4}

\newcommand\dn{\mathbb D_n}
\newcommand\dt{\mathbb D_3}

\newcommand\sm{\mathbb S_m}

\def\pf{\begin{proof}}
\def\epf{\end{proof}}

\theoremstyle{remark}

\begin{document}

\renewcommand{\baselinestretch}{1.2}

\thispagestyle{empty}

\title[pointed Hopf algebras]
{New techniques for pointed Hopf algebras}
\author[Andruskiewitsch and Fantino]{ Nicol\'as Andruskiewitsch and Fernando Fantino}
\thanks{This work was partially supported by
 ANPCyT-Foncyt, CONICET, Ministerio de Ciencia y Tecnolog\'\i a (C\'ordoba)  and Secyt (UNC)}
\address{\noindent Facultad de Matem\'atica, Astronom\'\i a y F\'\i sica,
Universidad Nacional de C\'ordoba. CIEM -- CONICET.
\newline \noindent Medina Allende s/n
(5000) Ciudad Universitaria, C\'ordoba, Argentina}
\email{andrus@famaf.unc.edu.ar} \email{fantino@famaf.unc.edu.ar}

\date{\today}

\begin{abstract}
We present techniques that allow to decide that the dimension of
some pointed Hopf algebras associated with non-abelian groups is
infinite. These results are consequences of \cite{AHS}. We
illustrate each technique with applications.
\end{abstract}
\maketitle

\centerline{\emph{Dedicado a Isabel Dotti y Roberto Miatello en su
sexag\'esimo cumplea\~nos}.}

\section*{Introduction}\label {0}

\subsection{}

Let $G$ be a finite group and let $\ydg$ be the category of
Yetter-Drinfeld modules over $\ku G$. The most delicate of the
questions raised by the Lifting Method for the classification of
finite-dimensional pointed Hopf algebras $H$ with $G(H)\simeq G$
\cite{AS1, AS-cambr}, is the following:

\begin{quote}\emph{Given $V\in \ydg$, decide when the Nichols algebra $\toba(V)$ is
finite-dimensional.} \end{quote}

Recall that a Yetter-Drinfeld module over the group algebra $\ku
G$ (or over $G$ for short) is a left $\ku G$-module and left $\ku
G$-comodule $M$ satisfying the compatibility condition
$\delta(g.m) = ghg^{-1} \otimes g.m$, for all $m\in M_h$, $g, h\in
G$. The list of all objects in $\ydg$ is known: any such is
completely reducible, and the class of irreducible ones is
parameterized by pairs $(\Oc, \rho)$, where $\Oc$ is a conjugacy
class in $G$ and $\rho$ is an irreducible representation of the
isotropy group $G^s$ of a fixed $s\in \Oc$. We denote the
corresponding Yetter-Drinfeld module by $M(\oc, \rho)$.

In fact, our present knowledge of Nichols algebras is still
preliminary. However, an important remark is that the Nichols
algebra $\toba(V)$ depends (as algebra and coalgebra) just on the
underlying braided vector space $(V,c)$-- see for example
\cite{AS-cambr}. This observation allows to go back and forth
between braided vector spaces and Yetter-Drinfeld modules. Indeed,
the same braided vector space could be realized as a
Yetter-Drinfeld module over different groups, and even in
different ways over the same group, or not at all. The braided
vector spaces that do appear as Yetter-Drinfeld modules over some
finite group are those coming from racks and 2-cocycles
\cite{AG1}. Thus, a comprehensive approach to the question above
would be to solve the following:

\begin{quote}\emph{Given a braided vector space $(V,c)$ determined by a rack
and a 2-cocycle, decide when $\dim \toba(V) < \infty$.}
\end{quote}

 But at the present moment and with the
exception of the diagonal case mentioned below, we know explicitly
very few Nichols algebras of braided vector spaces determined by
racks and 2-cocycles; see \cite{FK, MS, G1, AG1, G2}.

\subsection{}\label{small-lists} The braided vector spaces that appear as
Yetter-Drinfeld modules over some finite \emph{abelian} group are
the diagonal braided vector spaces. This leads to the following
question: \emph{Given a braided vector space $(V,c)$ of diagonal
type, decide when the Nichols algebra $\toba(V)$ is
finite-dimensional.} The full answer to this problem was given in
\cite{H-all}, see \cite{AS-adv, H-inv} for braided vector spaces
of Cartan type-- and \cite{AS6} for applications. These results on
Nichols algebras of braided vector spaces of diagonal type were in
turn used for more general pointed Hopf algebras. Let us fix a
non-abelian finite group $G$ and let $V\in \ydg$ irreducible. If
the underlying braided vector space contains a braided vector
subspace of diagonal type, whose Nichols algebra has infinite
dimension, then $\dim\toba(V)=\infty$. In turns out that, for
several finite groups considered  so far, many Nichols algebras of
irreducible Yetter-Drinfeld modules have infinite dimension; and
there are short lists of those not attainable by this method. See
\cite {G1, AZ, AF, AF2,FGV}.

\subsection{}\label{ahs} An approach of a different nature, inspired by
\cite{H-inv}, was  presented in \cite{AHS}. Let us consider  $V =
V_1 \oplus\dots\oplus V_\theta\in \ydg$, where the $V_i$'s are
irreducible. Then the Nichols algebra of $V$ is studied, under the
assumption that the $\toba(V_i)$ are known and finite-dimensional,
$1\le i \le \theta$. Under some circumstances, there is a Coxeter
group $\mathcal W$ attached to $V$, so that $\toba(V)$
finite-dimensional implies $\mathcal W$ finite. Although the
picture is not yet complete, the previous result implies that, for
a few $G$-- explicitly, $\st$, $\sk$, $\dn$-- the Nichols algebras
of some $V$ have infinite dimension. These applications rely on
the lists mentioned at the end of \ref{small-lists}.

\subsection{} The purpose of the present paper is to apply the
results in \ref{ahs} to discard more irreducible Yetter-Drinfeld
modules. Namely, let $V = V_1 \oplus V_2\in \ydga$, where $\Gamma
=\st$, $\sk$ or $\dn$, such that $\dim\toba(V)= \infty$ by
\cite[Section 4]{AHS}. Then there is a rack $(X, \trid)$ and a
cocycle $\q$ such that $(V,c) \simeq (\ku X, c_{\q})$. Let $G$ be
a finite group, let $\Oc$ be a conjugacy class in $G$, $s\in \Oc$,
$\rho\in \widehat{G^s}$ and $M(\Oc, \rho)\in \ydg$ the irreducible
Yetter-Drinfeld module corresponding to $(\Oc, \rho)$. We give
conditions on $(\Oc, \rho)$ such that $M(\Oc, \rho)$ contains a
braided vector subspace isomorphic to $(\ku X, c_{\q})$; thus,
necessarily, $\dim \toba(\Oc, \rho) = \infty$. We illustrate these
new techniques with several examples; see in particular Example
\ref{exa:12m} for one that can not be treated via abelian
subracks.

\subsection{} The facts glossed in the previous points strengthen
our determination to consider families of finite groups, in order
to discard those irreducible Yetter-Drinfeld modules over them
with infinite-dimensional Nichols algebra by the `subrack method'.
Natural candidates are the families of simple groups, or closely
related; cf. the classification of simple racks in \cite{AG1}. The
case of symmetric and alternating groups is treated in \cite{AZ,
AF, AF2, AFZ}; Mathieu groups in \cite{F}; other sporadic groups
in \cite{AFV}; some finite groups of Lie type with rank one in
\cite{FGV, FV}. Particularly, a list of only 9  types of pairs $(\Oc, \rho)$ for $\s_m$
whose Nichols algebras might be finite-dimensional is given in \cite{AFZ};
an analogous list of 7 pairs out of 1137 (for all 5 Mathieu simple groups) is given in \cite{F};
the sporadic groups $J_1$, $J_2$, $J_3$, $He$ and $Suz$ are shown to admit no
non-trivial pointed finite-dimensional Hopf algebra in \cite{AFV}.
Our new techniques are crucial for these results.

\subsection{} If for some finite group $G$ there is at most one irreducible
Yetter-Drinfeld module $V$ with finite-dimensional Nichols
algebra, then  \cite[Th. 4.2]{AHS} can be applied again. If the
conclusion is that $\dim \toba(V\oplus V) = \infty$, then we can
build a new rack together with a 2-cocycle realizing $V\oplus V$,
and investigate when a conjugacy class in another group $G'$
contains this rack, and so on.

\section{Notations and conventions}

The base field is $\mathbb C$ (the complex numbers).

\subsection{Braided vector spaces}\label{subsec:bvs}
A {\it braided vector space} is a pair $(V, c)$, where $V$ is a
vector space and $c: V\otimes V\to V\otimes V$ is a linear
isomorphism such that $c$ satisfies the braid equation: $
(c\otimes \id) (\id\otimes c) (c\otimes \id) = (\id\otimes c)
(c\otimes \id) (\id\otimes c)$.

\smallbreak Let $V$ be a vector space with a basis $(v_i)_{1\le
i\le \theta}$, let $(q_{ij})_{1\le i, j\le \theta}$ be a matrix of
non-zero scalars and let $c: V\otimes V\to V\otimes V$ be given by
$c(v_i\ot v_j) = q_{ij}v_j\ot v_i$. Then $(V,c)$ is a braided
vector space, called of \emph{diagonal type}.

\smallbreak Examples of braided vector spaces come from racks. A
\emph{rack} is a pair $(X,\trid)$ where $X$ is a non-empty set and
$\trid:X\times X\to X$ is a function-- called the multiplication,
such that $\phi_i:X\to X$, $\phi_i (j) := i\trid j$, is a
bijection for all $i\in X$, and
\begin{align}
\label{ccc2}& i\trid(j\trid k)=(i\trid j)\trid(i\trid k) \qquad
\text{ for all $i$, $j$, $k\in X$}.
\end{align}

For instance, a group $G$ is a rack with $x\trid y=x y x^{-1}$. In
this case, $j\trid i=i$ whenever $i\trid j=j$ and $i\trid i = i$
for all $i\in G$. We are mainly interested in subracks of $G$,
e.~g. in conjugacy classes in $G$.

Let $(X,\trid)$ be a rack. A function $\q: X\times X\to
\ku^{\times}$ is a \emph{2-cocycle} if $q_{i,j\trid k} \,
q_{j,k}=q_{i\trid j,i\trid k} \, q_{i,k}$, for all $i$, $j$, $k\in
X$. Then  $(\ku X, c_q)$ is a braided vector space, where $\ku X$
is the vector space with basis $e_k$, $k\in X$, and the braiding
is given by $$c_q(e_k\otimes e_l) = q_{k,l} \, e_{k\trid l}\otimes
e_k, \text{ for all }k, l\in X.$$ A subrack $T$ of $X$ is
\emph{abelian} if $k\trid l = l$ for all $k,l\in T$. If $T$ is an
abelian subrack of $X$, then $\ku T$ is a braided vector subspace
of $(\ku X, c_q)$ of diagonal type.

\begin{definition}\label{def:square-rack}
Let $X$ be a rack. Let $X_1$ and $X_2$ be two disjoint copies of
$X$, together with bijections $\varphi_i: X \to X_i$, $i=1,2$. The
\emph{square} of $X$ is the rack with underlying set the disjoint
union $X_1\coprod X_2$ and with rack multiplication
$$
\varphi_i(x) \trid \varphi_j(y) = \varphi_j(x \trid y),
$$
$x,y\in X$, $1\le i,j\le 2$. We denote the square of $X$ by
$X^{(2)}$. This is a particular case of an amalgamated sum of
racks, see e.~g. \cite{AG1}.
\end{definition}


\subsection{Yetter-Drinfeld modules}\label{subsec:yd}

We shall use the notation given in \cite{AF}. Let $G$ be a finite
group. We denote by $\vert g \vert$ the order of an element $g\in
G$; and by $\widehat{G}$ the set of isomorphism classes of
irreducible representations of $G$. We shall often denote a
representant of a class in $\widehat{G}$ with the same symbol as
the class itself.

Here is an explicit description of the irreducible Yetter-Drinfeld
module $M(\oc, \rho)$. Let $t_1 = s$, \dots, $t_{M}$ be a
numeration of $\oc$ and let $g_i\in G$ such that $g_i \trid s =
t_i$ for all $1\le i \le M$. Then  $ M(\oc, \rho) = \oplus_{1\le i
\le M}\, g_i\otimes V$, where $V$ is the vector space affording
the representation $\rho$. Let $g_iv := g_i\otimes v \in
M(\oc,\rho)$, $1\le i \le M$, $v\in V$. If $v\in V$ and $1\le i
\le M$, then the action of $g\in G$ is given by $g\cdot (g_iv) =
g_j(\gamma\cdot v)$, where $gg_i = g_j\gamma$, for some $1\le j
\le M$ and $\gamma\in G^s$, and the coaction is given by
$\delta(g_iv) = t_i\otimes g_iv$. Then $M(\oc, \rho)$ is a braided
vector space with braiding $ c(g_iv\otimes g_jw) = g_h(\gamma\cdot
w)\otimes g_iv$, for any $1\le i,j\le M$, $v,w\in V$, where
$t_ig_j = g_h\gamma$ for unique $h$, $1\le h \le M$ and $\gamma
\in G^s$. Since $s\in Z(G^s)$, the center of $G^{s}$, the Schur
Lemma implies that
\begin{equation}\label{schur} s \text{ acts by a scalar $q_{ss}$
on } V.
\end{equation}

\begin{lema}\label{trivialbraiding}
If $U$ is a subspace of $W$ such that $c(U\otimes U) = U\otimes U$
and $\dim \toba(U) =\infty$, then $\dim \toba(W) =\infty$.\qed
\end{lema}

\begin{lema}\label{odd}\cite[Lemma 2.2]{AZ}  Assume that $s$ is \emph{real} (i.~e. $s^{-1}\in \oc$).
If $\dim\toba(\Oc, \rho)< \infty$, then $q_{ss} = -1$ and $s$ has
even order.\qed
\end{lema}

Let $\trasp\in \sm$ be a product of $n_j$ disjoint cycles of
length $j$, $1\le j \le m$. Then the type of $\trasp$ is the
symbol $(1^{n_1}, 2^{n_2}, \dots, m^{n_m})$. We may omit $j^{n_j}$
when $n_j = 0$. The conjugacy class $\Oc_\trasp$ of $\trasp$
coincides with the set of all permutations in $\sm$ with the same
type as $\trasp$; we may use the type as a subscript of a
conjugacy class as well. If some emphasis is needed, we add a
superscript $m$ to indicate that we are taking conjugacy classes
in $\sm$, like $\Oc_j^m$ for the conjugacy class of $j$-cycles in
$\sm$.

\section{A technique from the dihedral group $\dn$, $n$ odd}

Let $n>1$ be an odd integer. Let $\dn$ be the dihedral group of
order $2n$, generated by $x$ and $y$ with defining relations $x^2
= e = y^n$ and $xyx = y^{-1}$. Let $\Oc_x$ be the conjugacy class
of $x$ and let $\sgn \in \widehat{\dn^x}$ be the sign
representation ($\dn^x=\langle x \rangle \simeq \mathbb Z_2$). The
goal of this Section is to apply the next result, cf. \cite[Th.
4.8]{AHS}, or \cite[Th. 4.5]{AHS} for $n =3$.

\begin{theorem}\label{theorem:dn-ahs}
The Nichols algebra  $\toba(M(\Oc_{x}, \sgn) \oplus M(\Oc_{x},
\sgn))$ has infinite dimension. \qed
\end{theorem}

Note that $M(\Oc_{x}, \sgn) \oplus M(\Oc_{x}, \sgn)$ is isomorphic
as a braided vector space to $(\ku X_n, \q)$, where

\begin{itemize}
    \item $X_n$ is the rack  with $2n$ elements $s_i$, $t_j$, $i,j\in
    \Z/n$, and with structure
    $$
s_i \trid s_j = s_{2i - j}, \quad s_i \trid t_j = t_{2i - j},
\quad t_i \trid s_j = s_{2i - j},\quad t_i \trid t_j = t_{2i - j},
\quad i,j\in \Z/n;
    $$
    \item $\q$ is the constant cocycle $\q \equiv -1$.
\end{itemize}

\medbreak If $d$ divides $n$, then $X_d$ can be identified with a
subrack of $X_n$. Hence, it is enough to consider braided vector
spaces $(\ku X_p, \q)$, with $p$ an odd prime.

\medbreak We fix a finite group $G$ with the rack structure given
by conjugation $x \trid y=xyx^{-1}$, $x$, $y\in G$. Let $\oc$ be a
conjugacy class in $G$.

\begin{definition}\label{defi:dp} Let $p>1$ be an integer. A family $(\mu_i)_{i\in
\Z/p}$ of distinct elements of $G$ is \emph{of type $\D_p$} if
\begin{equation}\label{eqn:dp}
\mu_i \trid \mu_j = \mu_{2i - j}, \quad i,j\in \Z/p.
\end{equation}
Let $(\mu_i)_{i\in \Z/p}$ and $(\nu_i)_{i\in \Z/p}$ be two
families of type $\D_p$ in $G$, such that  $\mu_i\neq\nu_j$ for
all $i,j\in \Z/p$. Then $(\mu, \nu) := (\mu_i)_{i\in \Z/p}\cup
(\nu_i)_{i\in \Z/p}$ is \emph{of type $\D_p^{(2)}$} if
\begin{equation}\label{eqn:dp2}
\mu_i \trid \nu_j = \nu_{2i - j}, \quad\nu_i \trid \mu_j = \mu_{2i
- j}, \quad i,j\in \Z/p.
\end{equation}

\end{definition}

\bigbreak

It is useful to denote $i\trid j = 2i - j$, for $i,j\in \Z/p$.


We state some consequences of this definition for further use.

\begin{obs} If $(\mu_i)_{i\in \Z/p}$ is of type $\D_p$ then
\begin{align}\label{eq:mu^menosuno}
\mu_i^{-1} \trid \mu_j &= \mu_{2i - j}, &\qquad \mu_i \trid
\mu_j^{-1} &= \mu_{2i - j}^{-1}, &\qquad \mu_i^{-1} \trid
\mu_j^{-1} &= \mu_{2i - j}^{-1},
\\\label{eq:mu^r}
\mu_i^{k} \trid \mu_j &= \mu_{2i - j}, &\qquad \mu_i \trid
\mu_j^{k} &= \mu_{2i - j}^{k}, &\qquad \mu_i^{k} \trid \mu_j^{k}
&= \mu_{2i - j}^{k},
\end{align}
for all $i,j\in \Z/p$, and for all $k$ odd.
\end{obs}

\begin{obs} Assume that $p$ is odd.
If $(\mu, \nu)= (\mu_i)_{i\in \Z/p}\cup (\nu_i)_{i\in \Z/p}$ is of
type $\D_p^{(2)}$, then for all $i$, $j$,
\begin{align}\label{eq:mu_i^2=nu_j^2}
\mu_i^2=\mu_j^2,\quad \nu_i^2=\nu_j^2,\quad  \mu_i^2
\nu_j=\nu_j\mu_i^2, \quad \nu_i^2 \mu_j=\mu_j\nu_i^2.
\end{align}
Indeed, $\mu_h^2\mu_j = \mu_j\mu_h^2$, hence $\mu_{2h - j}^2 =
\mu_h\mu_j^2 \mu_h^{-1} = \mu_j^2$. Take now $h = \dfrac{i+j}2$.
\end{obs}

\begin{lema}\label{le:dp2}
If $(\mu, \nu)= (\mu_i)_{i\in \Z/p}\cup (\nu_i)_{i\in \Z/p}$ is of
type $\D_p^{(2)}$, then
\begin{itemize}
\item[(i)] $\mu_k\mu_l=\mu_{t(l-k)+k} \,\, \mu_{t(l-k)+l}$,
\item[(ii)] $\mu_k\nu_l=\mu_{2t(l-k)+k} \,\, \nu_{2t(l-k)+l}$,
\item[(iii)] $\mu_k\nu_l=\nu_{(2t+1)(l-k)+k} \,\,
\mu_{(2t+1)(l-k)+l}$,
\end{itemize}
for all $k$, $l$, $t\in \Z/p$.
\end{lema}

Notice that we have the analogous relations interchanging $\mu$ by
$\nu$.

\pf We proceed by induction on $t$. We will prove (i); (ii) and
(iii) are similar. The result is obvious when $t=0$. Since
$\mu_k\mu_l=\mu_l \,\,\mu_{l\trid k}$, then the result holds for
$t=1$. Let us suppose that (i) holds for every $s\leq t$. Now,
\begin{align*}
\mu_k\mu_l&=\mu_{t(l-k)+k}\,\,\mu_{t(l-k)+l}
\\
&= \mu_{t(l-k)+l}\,\,\mu_{(t(l-k)+l)\trid (t(l-k)+k)} =
\mu_{(t+1)(l-k)+k}\,\,\mu_{(t+1)(l-k)+l}
\end{align*}
by the recursive hypothesis. \epf

\begin{lema}\label{le:dp2:nu:mu} Assume that $p$ is odd.
If $(\mu, \nu)$ is of type $\D_p^{(2)}$, then for $i\in \Z/p$,
\begin{align}
\label{eq:mu:nu} \mu_i \nu_i&= \mu_0 \nu_0,\\
\label{eq:nu:mu} \nu_i \mu_i&= \nu_0 \mu_0.
\end{align}
\end{lema}

\pf Let $i$, $j\in \Z/p$, with $i\neq j$. If we write (ii) of
Lemma \ref{le:dp2} with $k=i$, $l=j$ and $t=-1/2$ we have that
$\mu_i\nu_j=\mu_{2i-j}\nu_{i}$. Thus, $\mu_i\nu_i \nu_j^2=\mu_i
\nu_j \nu_j \nu_i=\mu_{2i-j}\nu_{i} \nu_i \nu_{2i-j}=\mu_{2i-j}
\nu_{2i-j}\nu_{i}^2$, and, by \eqref{eq:mu_i^2=nu_j^2},
\begin{align*}
\mu_i\nu_i=\mu_{2i-j} \nu_{2i-j}.
\end{align*}
Now \eqref{eq:mu:nu} follows taking $j=2i$. Now
\eqref{eq:nu:mu} follows from \eqref{eq:mu:nu} by \eqref{eqn:dp2}. \epf

\smallbreak We now set up some notation that will be used in the
rest of this section. Let $(\mu_i)_{i\in \Z/p}$ be a family of
type $\D_p$ in $G$, with $p$ odd. Set
\begin{align}
\label{eqn:dp-gs} g_i &= \mu_{i/2},
\\\label{eqn:dp-as} \alpha_{ij} &= g_{i\trid
j}^{-1}\,\mu_i\, g_j = \mu_{i - j/2}^{-1} \, \mu_{i} \, \mu_{j/2},
\end{align}
for all $i,j\in \Z/p$. Then $$ g_i\trid \mu_0 = \mu_i, \qquad
\alpha_{ij}\in G^{\,\mu_0}, \qquad i,j\in \Z/p.
$$

Let now $(\mu,\nu)$ be of type $\D_p^{(2)}$ and suppose that there
exists $\gi\in G$ such that $\gi\trid \mu_0 =\nu_0$. Set
\begin{align}
\label{eqn:dp-fs}f_i &= \nu_{i/2}\,\gi,
\\\label{eqn:dp-bs}
\beta_{ij} &= f_{i\trid j}^{-1}\,\mu_i \,f_j = \gi^{-1}\,\nu_{i -
j/2}^{-1}\, \mu_{i}\, \nu_{j/2}\,\gi,
\\\label{eqn:dp-gas}
\gamma_{ij} &= g_{i\trid j}^{-1}\,\nu_i \,g_j = \mu_{i - j/2}^{-1}
\,\nu_{i} \,\mu_{j/2},
\\\label{eqn:dp-ds}
\delta_{ij} &= f_{i\trid j}^{-1}\,\nu_i \,f_j = \gi^{-1}\,\nu_{i -
j/2}^{-1} \,\nu_{i} \,\nu_{j/2}\,\gi.
\end{align}
Then
$$f_i\trid \mu_0 = \nu_i, \qquad \beta_{ij},\, \gamma_{ij},\,
\delta_{ij} \in G^{\,\mu_0}, \qquad i,j\in \Z/p.$$

\medbreak \emph{We assume from now on that $p$ is an odd prime.}
This is required in the proof of the next lemma, needed for the
main result of the section.

\begin{lema}\label{le:dp2:a:b:g:d}
Let $(\mu, \nu)= (\mu_i)_{i\in \Z/p}\cup (\nu_i)_{i\in \Z/p}$ be
of type $\D_p^{(2)}$, and suppose that there exists $\gi\in G$
such that $\gi\trid \mu_0 =\nu_0$. Let $g_i$ and $f_i$ be as in
\eqref{eqn:dp-gs} and \eqref{eqn:dp-fs}, respectively. Then, for
all $i,j\in \Z/p$,
\begin{itemize}
\item[(a)] $\alpha_{ij}=\delta_{ij}=\mu_0$,
\item[(b)] $\beta_{ij} = \gi^{-1}\mu_0 \gi$,
\item[(c)] $\gamma_{ij} = \nu_0$.
\end{itemize}
\end{lema}

\pf Let $k$, $l$ be in $\Z/p$. Then, for all $r\in \Z/p$, we have
\begin{equation}
\label{eq:mu:nu:heck} \mu_k\mu_l = \mu_{k+r}\mu_{l+r}, \quad
\mu_k\nu_l = \mu_{k+r}\nu_{l+r}, \quad\mu_k\nu_l =
\nu_{k+r}\mu_{l+r}.
\end{equation}
This follows from \eqref{eq:mu_i^2=nu_j^2} and Lemma
\ref{le:dp2:nu:mu} (when $k=l$), and Lemma \ref{le:dp2} (when
$k\neq l$). There are similar equalities interchanging $\mu$'s and
$\nu$'s. Now
\begin{align*}
\alpha_{ij} &= \mu_{i - j/2}^{-1}\,\mu_{i}\, \mu_{j/2}
\overset{\eqref{eq:mu:nu:heck}}=\mu_0, \\ \delta_{ij} &=
\gi^{-1}\nu_{i - j/2}^{-1}\,\nu_{i}\, \nu_{j/2}\,\gi
\overset{\eqref{eq:mu:nu:heck}}=\gi^{-1}\nu_0\,\gi=\mu_0,
\\
\beta_{ij} &=\gi^{-1}  \nu_{i - j/2}^{-1}\,\mu_{i}\, \nu_{j/2}\,
\gi \overset{\eqref{eq:mu:nu:heck}}=\gi^{-1}\,\mu_0\,\gi,
\\
\gamma_{ij}&=\mu_{i - j/2}^{-1}\,\nu_{i}\, \mu_{j/2}
\overset{\eqref{eq:mu:nu:heck}}= \mu_{i - j/2}^{-1}\,\mu_{i - j/2}
\nu_{0} =\nu_0,
\end{align*}
and the Lemma is proved. \epf

\medbreak

We can now prove one of the main results of this paper.

\begin{theorem}\label{theorem:dp-cor}
Let $(\mu, \nu) = (\mu_i)_{i\in \Z/p}\cup (\nu_i)_{i\in \Z/p}$ be
a family of elements in $G$ with $\mu_0\in \oc$. Let $(\rho, V)$
be an irreducible representation of the centralizer $G^{\,\mu_0}$.
We assume that
\begin{itemize}
\item[(H1)] $(\mu,\nu)$ is of type $\D_p^{(2)}$;
\item[(H2)] $(\mu,\nu)\subseteq \oc$, with $\gi\in G$ such
that $\gi\trid \mu_0 =\nu_0$;
\item[(H3)] $q_{\mu_0\mu_0}=-1$;
\item[(H4)] there exist $v,w\in V - 0$ such that,
\begin{align}
\label{eqn:dp-b-act} \rho(\gi^{-1}\mu_0 \gi) w &= -w,
\\\label{eqn:dp-ga-act} \rho(\nu_0) v &= -v.
\end{align}
\end{itemize} Then $\dim\toba(\oc,\rho)=\infty$.
\end{theorem}

\pf We keep the notation \eqref{eqn:dp-fs}--\eqref{eqn:dp-ds}
above. Let $v, w\in V - 0$ as in (H4) and let $W:=$ span-$\{g_iv:
i\in \Z/p\}\cup \{f_iw: i\in \Z/p\}$. Let $\Psi: \ku X_p \to W$ be
given by $\Psi(s_i) = g_iv$, $\Psi(t_i) = f_iw$, $i\in \Z/p$.
Since the elements $\mu_i$ and $\nu_j$ are all different, $\Psi$
is a linear isomorphism. We claim that $W$ is a braided vector
subspace of $M(\Oc, \rho)$ and that $\Psi$ is an isomorphism of
braided vector spaces. We compute the braiding in $W$:
\begin{align*}
c(g_iv\ot g_jv)& = \mu_ig_jv \ot g_iv = g_{i\trid j}\alpha_{ij}v
\ot g_iv \overset{\text{(H3)}}= -
g_{i\trid j} v \ot g_iv,\\
c(g_iv\ot f_jw)& = \mu_if_jw \ot g_iv = f_{i\trid j}\beta_{ij}w
\ot g_iv \overset{\eqref{eqn:dp-b-act}}= -
f_{i\trid j} w \ot g_iv,\\
c(f_iw\ot g_jv)& = \nu_ig_jv \ot f_iw = g_{i\trid j}\gamma_{ij}v
\ot f_iw \overset{\eqref{eqn:dp-ga-act}}= -
g_{i\trid j} v \ot f_iw,\\
c(f_iw\ot f_jw)& = \nu_if_jw \ot f_iw = f_{i\trid j}\delta_{ij}w
\ot f_iw \overset{\text{(H3)}}= - f_{i\trid j} w \ot f_iw,
\end{align*}
by Lemma \ref{le:dp2:a:b:g:d}. The claim is proved. Hence,
$\dim\toba(W)=\infty$ by Theorem \ref{theorem:dn-ahs}. Now the
Theorem follows from Lemma \ref{trivialbraiding}. \epf

As a consequence of Theorem \ref{theorem:dp-cor}, we can state a
very useful criterion.

\begin{cor}\label{coro:dp-cor}
Let $G$ be a finite group, $\mu_i$, $0 \leq i \leq p-1$, distinct
elements in $G$, with $p$ an odd prime. Let us suppose that there
exists $k\in \Z$ such that $\mu_0^k\neq\mu_0$ and $\mu_0^k\in
\oc$, the conjugacy class of $\mu_0$. Let $\rho=(\rho,V)\in
\widehat{G^{\,\mu_0}}$. Assume further that
\begin{itemize}
\item[(i)] $(\mu_i)_{i\in \Z/p}$ is of type $\D_p$,
\item[(ii)] $q_{\mu_0\mu_0}=-1$.
\end{itemize}
Then $\dim\toba(\oc,\rho)=\infty$.
\end{cor}

\pf We may assume that $1< k < |\mu_0|$. By hypothesis (ii),  the
order of $\mu_0$ is even; hence $k$ is odd, say $k=2t+1$, with
$t\geq 1$. Let $\nu_i:=\mu_i^{k}$, $0\leq i \leq p-1$, and pick
$\gi\in G$ such that $\gi\trid \mu_0 =\mu_0^{k}$. Set $(\mu, \nu)=
(\mu_i)_{i\in \Z/p}\cup (\nu_i)_{i\in \Z/p}$. Clearly
$(\mu,\nu)\subseteq \oc$. We claim that $(\mu,\nu)$ is of type
$\D_p^{(2)}$. Indeed, using (i) it is easy to see that
$(\mu_i)_{i\in \Z/p}\cup (\nu_i)_{i\in \Z/p}$ are all distinct.
Then the claim follows by \eqref{eq:mu^r}.

\smallbreak It remains to check the hypothesis (H4) of Theorem
\ref{theorem:dp-cor}. As $\gi\mu_0 \gi^{-1}=\mu_0^k$,
$\gi^{l}\mu_0 \gi^{-l}=\mu_0^{k^l}$, for all $l\geq 0$. In
particular, $$\gi^{-1}\mu_0 \gi=\gi^{|\gi|-1}\mu_0
\gi^{-|\gi|+1}=\mu_0^{k^{|\gi|-1}}.$$ Then, since
$q_{\mu_0\mu_0}=-1$ and $k$ is odd, we see that
$\rho(\gi^{-1}\mu_0 \gi)=-\Id$. Hence \eqref{eqn:dp-b-act} holds,
for any $w\in V-0$. Also,
$\rho(\nu_0)=\rho(\mu_0^k)=(-\Id)^k=-\Id$, because $k$ is odd;
thus, \eqref{eqn:dp-ga-act} holds for any $v\in V-0$. Thus, for
any  $v$, $w$ in $V-0$, we are in the conditions of Theorem
\ref{theorem:dp-cor}. Then $\dim\toba(\oc,\rho)=\infty$.\epf

\begin{example}\label{co:type6k^nj}
Let $m\geq 6$. Let $\sigma\in \sm$ of type
$(1^{n_1},2^{n_2},\dots,m^{n_m})$, $\oc$ the conjugacy class of
$\sigma$ and $\rho \in \widehat{\s_m^{\sigma}}$. If there exists
$j$, $1\leq j \leq m$, such that
\begin{itemize}
    \item $2p$ divides $j$, for some odd prime $p$, and
    \item $n_j\geq 1$;
\end{itemize}
then $\dim \toba(\oc,\rho)=\infty$.
\end{example}

Before proving the Example, we state a more general Lemma that
might be of independent interest. Here $p$ is no longer an odd
prime.

\begin{lema}\label{claim:type6k^nj} Let $m, p\in \Z_{>1}$.  Let $\sigma\in \sm$ of type
$(1^{n_1},2^{n_2},\dots,m^{n_m})$ and $\oc$ the conjugacy class of
$\sigma$. If there exists $j\neq 4$, $1\leq j \leq m$, such that
\begin{itemize}
    \item $2p$ divides $j$,  and
    \item $n_j\geq 1$;
\end{itemize}
then $\oc$ contains a subrack of type $\D_p^{(2)}$.
\end{lema}

\pf Let $j=2p\,\kappa$, with $\kappa\geq 1$. Let $\alpha=(i_{1} \,
i_{2} \, \cdots \, i_{j})$ be a $j$-cycle that appears in the
decomposition of $\sigma$ as product of disjoint cycles and define
\begin{align*}
{\bf I}:=(i_{1}\, i_{3}\, i_{5} \, \cdots \, i_{j-1})\quad \text{
and } \quad {\bf P}:=(i_{2}\, i_{4}\, i_{6} \, \cdots \, i_{j}).
\end{align*}
We claim that
\begin{itemize}
\item[(a)] ${\bf I}$ and ${\bf P}$ are disjoint $p\kappa$-cycles,
\item[(b)] $\alpha^2={\bf I} {\bf P}$,
\item[(c)] $\alpha  {\bf I} \alpha^{-1}={\bf P}$, (and then $\sigma  {\bf I} \sigma^{-1}={\bf
P}$),
\item[(d)] ${\bf P}^t \alpha {\bf P}^{t}=\alpha^{2t+1}$,
${\bf P}^t \alpha^{-1} {\bf P}^{t}=\alpha^{2t-1}$, for all
integers $t$.
\end{itemize}
The first two items are clear, while (c) follows from the
well-known formula $\alpha(l_1\, l_2\,\dots\, l_k) \alpha^{-1} =
(\alpha(l_1)\, \alpha(l_2)\,\dots\, \alpha(l_k))$. (d). By (c),
${\bf P}^t=\alpha  {\bf I}^t \alpha^{-1}$. Then ${\bf P}^t \alpha
{\bf P}^t=\alpha {\bf I}^t {\bf P}^t$; by (b), ${\bf P}^t \alpha
{\bf P}^t= \alpha \alpha^{2t}$, as claimed.

We define
\begin{align}
\sigma_i:={\bf P}^{i\kappa} \sigma {\bf P}^{-i\kappa}, \qquad
0\leq i \leq p-1.
\end{align}
Notice that $\sigma_i={\bf P}^{i\kappa} \alpha {\bf
P}^{-i\kappa}\,\,\widetilde{\sigma}$, where
$\widetilde{\sigma}:=\alpha^{-1}\sigma$. The elements
$(\sigma_i)_{i\in \Z/p}$ are all distinct; indeed, if
$\sigma_i=\sigma_l$, with $i$, $l\in \Z/p$, then ${\bf
P}^{i\kappa} \sigma {\bf P}^{-i\kappa}={\bf P}^{l\kappa} \sigma
{\bf P}^{-l\kappa}$, i.~e. ${\bf P}^{(i-l)\kappa} \sigma {\bf
P}^{-(i-l)\kappa}=\sigma$, which implies that
$i_2=\sigma(i_1)={\bf P}^{(i-l)\kappa} \sigma {\bf
P}^{-(i-l)\kappa}(i_1)={\bf P}^{(i-l)\kappa}
(i_2)=i_{2(i-l)\kappa+2}$, and this means that $2(i-l)\kappa=0$ in
$\Z/j$. Thus $i=l$, as desired.

\medbreak We claim that $(\sigma_i)_{i\in \Z/p}$ is of type
$\D_p$. If $i$, $l\in \Z/p$, then
\begin{align*}
\sigma_i\trid \sigma_l&= {\bf P}^{i\kappa} \sigma {\bf
P}^{-i\kappa} \,\, {\bf P}^{l\kappa} \sigma {\bf P}^{-l\kappa}
\,\, {\bf P}^{i\kappa} \sigma^{-1}
{\bf P}^{-i\kappa}\\
&={\bf P}^{i\kappa}\, \alpha\, {\bf P}^{-i\kappa} \,\, {\bf
P}^{l\kappa} \,\alpha\, {\bf P}^{-l\kappa} \,\, {\bf P}^{i\kappa}
\alpha^{-1} {\bf
P}^{-i\kappa} \,\,\widetilde{\sigma}\\
&= {\bf P}^{(2i-l)\kappa}\, {\bf P}^{(l-i)\kappa}\, \alpha\, {\bf
P}^{(l-i)\kappa} \, \alpha \,{\bf P}^{(i-l)\kappa} \,
 \alpha^{-1}\, {\bf P}^{(i-l)\kappa}\, {\bf P}^{-(2i-l)\kappa}\,\,\widetilde{\sigma}\\
&= {\bf P}^{(2i-l)\kappa}\, \alpha^{2(l-i)\kappa+1} \, \alpha \,
\alpha^{2
(i-l)\kappa -1} \, {\bf P}^{-(2i-l)\kappa}\,\,\widetilde{\sigma}\\
&= {\bf P}^{(2i-l)\kappa}\, \alpha \, {\bf
P}^{-(2i-l)\kappa}\,\,\widetilde{\sigma}= {\bf P}^{(2i-l)\kappa}\,
\sigma \, {\bf P}^{-(2i-l)\kappa}=\sigma_{i\trid l},
\end{align*}
by (d), and the claim follows. Finally, the family of type
$\D_p^{(2)}$ we are looking for is $(\sigma_i)_{i\in \Z/p}\cup
(\sigma_i^{-1})_{i\in \Z/p}$. It remains to show that
$\sigma_t\neq\sigma_l^{-1}$ for all $t$,
$l\in\Z_p $. If  $\sigma_t = \sigma_l^{-1}$, then $\sigma_t^2(i_1) = \sigma_l^{-2}(i_1)$,
that is $i_3 = i_{j-1}$, a contradiction to the hypothesis $j\neq4$.
\epf

\emph{Proof of the Example \ref{co:type6k^nj}.}  We may assume
that $q_{\sigma\sigma}= -1$, by Lemma \ref{odd}. By Lemma
\ref{claim:type6k^nj}, we have a family $(\sigma_i)_{i\in \Z/p}$
of type $\D_p$, with $\sigma_0 = \sigma$. Now Corollary
\ref{coro:dp-cor} applies, with $\mu_0 = \sigma_0$, $k = \vert
\sigma_0\vert -1$. Thus $\dim \toba(\oc,\rho)=\infty$. \qed

\section{A technique from the symmetric group
$\st$}\label{nuevatecnica:R2} We study separately the case $p=3$
because of the many applications found. In this setting,
$\dt\simeq \st$ and $\oc_x = \oc_2^3 =\{(1\, 2),\, (2\, 3),\, (1\,
3)\}$ is the conjugacy class of transpositions in $\st$. The rack
$X_3$ is described as  a set of 6 elements $X_3 = \{x_1, x_2, x_3,
y_1, y_2, y_3\}$ with the multiplication
\begin{align*} x_i \trid x_j=x_k, \quad y_i \trid y_j=y_k, \quad
x_i \trid y_j=y_k, \quad y_i \trid x_j=x_k,
\end{align*}
for $i$, $j$, $k$, all distinct or all equal.

\subsection{Families of type $\D_3$ and
$\D_3^{(2)}$}\label{subsubsec:r2}

We fix a finite group $G$ and $\oc$ a conjugacy class in $G$. Our
aim is to give criteria to detect when $\oc$ contains a subrack
isomorphic to $X_3$.

\begin{definition} Let $\sigma_1$, $\sigma_2$,
$\sigma_3\in G$ distinct. We say  that $(\sigma_i)_{1\leq i \leq
3}$ is of type $\D_3$ if
\begin{align}\label{R2} \sigma_i\trid \sigma_j=\sigma_k,
\qquad \text{ where $i$, $j$, $k$ are all distinct}.
\end{align}
\end{definition}

The requirement \eqref{R2} consists of 6 identities, but actually
3 are enough.

\begin{remark}\label{obs:prop1}
If
\begin{align}
\label{123} \sigma_1\trid \sigma_2&=\sigma_3,\\
\label{132} \sigma_1\trid \sigma_3&=\sigma_2,\\
\label{231} \sigma_2\trid \sigma_3&=\sigma_1,
\end{align}
then $(\sigma_i)_{1\leq i \leq 3}$ is of type $\D_3$. \qed
\end{remark}

Here is a characterization of $\D_3$ families.

\begin{prop}\label{prop:caracR2}
Let $\sigma_1$, $\sigma_2\in \oc$. Define $\sigma_3:=\sigma_1
\trid \sigma_2$. Then $(\sigma_i)_{1\leq i \leq 3}$ is of type
$\D_3$ if and only if \begin{align}
\label{eq1:prop1} &\text{$\sigma_1\not\in G^{\sigma_2}$},\\
\label{eq2:prop1} &\text{$\sigma_1^2\in G^{\sigma_2}$},\\
\label{eq3:prop1} &\text{$\sigma_1=\sigma_2\trid(\sigma_1 \trid
\sigma_2)$}.
\end{align}
\end{prop}

\pf The definition of $\sigma_3$ is equivalent to \eqref{123} and
\eqref{eq3:prop1} is equivalent to \eqref{231}. Assume that
$(\sigma_i)_{1\leq i \leq 3}$ is of type $\D_3$. As $\sigma_3\neq
\sigma_2$, $\sigma_1\not\in G^{\sigma_2}$. Also, $\sigma_1^2 \trid
\sigma_2=\sigma_1\trid (\sigma_1 \trid \sigma_2)=\sigma_1 \trid
\sigma_3=\sigma_2$. Hence $\sigma_1^2\in G^{\sigma_2}$.

Conversely, if $\sigma_1\not\in G^{\sigma_2}$, then $\sigma_1\neq
\sigma_2$, $\sigma_2\neq \sigma_3$. From \eqref{eq1:prop1} and
\eqref{eq3:prop1}, we see that $\sigma_1\neq \sigma_3$. It remains
to check \eqref{132}: $\sigma_1\trid \sigma_3=\sigma_1^2\trid
\sigma_2=\sigma_2$. \epf

\begin{definition}
Let $\sigma_1$, $\sigma_2$, $\sigma_3$, $\tau_1$, $\tau_2$,
$\tau_3\in G$ be distinct elements. We say  that
$(\sigma,\tau)=(\sigma_1,\sigma_2,\sigma_3,\tau_1,\tau_2, \tau_3)$
is of type $\D_3^{(2)}$, if $(\sigma_i)_{1\leq i \leq 3}$ and
$(\tau_j)_{1\leq j \leq 3}$ are of type $\D_3$, and
\begin{align}\label{R2^2}
\sigma_i\trid \tau_j=\tau_k, \qquad  \tau_i\trid
\sigma_j=\sigma_k,
\end{align}
where $i$, $j$, $k$ are either all equal, or all distinct.
\end{definition}

The requirement \eqref{R2^2} consists of 18 identities, but less
are enough. We begin by a first reduction.

\begin{lema}\label{lema:st}
Let $(\sigma_i)_{1\leq i \leq 3}$ and $(\tau_j)_{1\leq j \leq 3}$
such that \eqref{123}, \eqref{132}, \eqref{231} hold for $\sigma$
and for $\tau$. If
\begin{align}
\label{st111} \sigma_1\trid \tau_1&=\tau_1,\\
\label{st123} \sigma_1\trid \tau_2&=\tau_3,\\
\label{st213} \sigma_2\trid \tau_1&=\tau_3,
\end{align}
also hold, then $\sigma_i\trid \tau_i=\tau_i$, $1\leq i \leq 3$,
and $\sigma_i\trid \tau_j=\tau_k$, for all $i$, $j$, $k$ distinct.
\end{lema}

\pf We have to prove
\begin{align}
\label{st132} \sigma_1\trid \tau_3&=\tau_2,\\
\label{st333} \sigma_3\trid \tau_3&=\tau_3,\\
\label{st222} \sigma_2\trid \tau_2&=\tau_2,\\
\label{st312} \sigma_3\trid \tau_1&=\tau_2,\\
\label{st321} \sigma_3\trid \tau_2&=\tau_1,\\
\label{st231} \sigma_2\trid \tau_3&=\tau_1,
\end{align}

The identity \eqref{st132} holds because $ \sigma_1\trid \tau_3=
\sigma_1 \trid(\tau_1 \trid \tau_2)=\tau_1 \trid \tau_3=\tau_2$;
in turn, \eqref{st333} and \eqref{st222} hold because
\begin{align*}
\sigma_3\trid \tau_3= (\sigma_2 \trid \sigma_1) \trid (\sigma_2
\trid \tau_1)=\sigma_2 \trid( \sigma_1 \trid \tau_1)=\sigma_2
\trid \tau_1=\tau_3,\\
\sigma_2\trid \tau_2= (\sigma_1 \trid \sigma_3) \trid (\sigma_1
\trid \tau_3)=\sigma_1 \trid( \sigma_3 \trid \tau_3)=\sigma_1
\trid \tau_3=\tau_2.
\end{align*}
Also, $\sigma_3\trid \tau_1= (\sigma_1 \trid \sigma_2) \trid
(\sigma_1 \trid \tau_1)=\sigma_1 \trid( \sigma_2 \trid
\tau_1)=\sigma_1 \trid \tau_3=\tau_2$, showing \eqref{st312}.
Finally, $\sigma_3\trid \tau_2 = \sigma_3 \trid( \sigma_1 \trid
\tau_3)=\sigma_2 \trid( \sigma_3 \trid \tau_3)=\sigma_2 \trid
\tau_3=\sigma_2 \trid( \tau_1 \trid \tau_2)=\tau_3 \trid
\tau_2=\tau_1$, proving \eqref{st321} and \eqref{st231}. \epf

Therefore, given 6 distinct elements $\sigma_1$, $\sigma_2$, $\sigma_3$,
$\tau_1$, $\tau_2$, $\tau_3\in G$, if the 12 identities:
\eqref{123}, \eqref{132}, \eqref{231}, for $\sigma$ and for
$\tau$, \eqref{st111}, \eqref{st123}, \eqref{st213}, and the
analogous identities
\begin{align}
\label{ts111} \tau_1\trid \sigma_1&=\sigma_1,\\
\label{ts123} \tau_1\trid \sigma_2&=\sigma_3,\\
\label{ts213} \tau_2\trid \sigma_1&=\sigma_3,
\end{align}
hold, then $(\sigma, \tau)$ is of type $\D_3^{(2)}$. But we can
get rid of 3 of these 12 identities.

\begin{prop}
Let $\sigma_1$, $\sigma_2$, $\sigma_3$, $\tau_1$, $\tau_2$,
$\tau_3\in G$, all distinct, such that \eqref{123}, \eqref{132}, \eqref{231},
hold for $\sigma$ and for $\tau$, as well as the identities
\eqref{st111}, \eqref{st213} and \eqref{ts123}. Then $(\sigma,
\tau)$ is of type $\D_3^{(2)}$.
\end{prop}

\pf By Lemma \ref{lema:st}, it is enough to check \eqref{st123},
\eqref{ts111} and \eqref{ts213}. First, \eqref{ts111} holds
because $\tau_1=\sigma_1 \trid \tau_1=\sigma_1 \tau_1
\sigma_1^{-1}$. If $\tau_1$ acts on both sides of \eqref{st213},
then $\tau_2=\tau_1\trid \tau_3=(\tau_1 \trid \sigma_2) \trid
(\tau_1 \trid \tau_1)=\sigma_3 \trid \tau_1$; if now $\sigma_1$
acts on the last, then
\begin{align*} \sigma_1 \trid \tau_2=(\sigma_1 \trid \sigma_3)
\trid (\sigma_1 \trid \tau_1)=\sigma_2 \trid \tau_1
\overset{\eqref{st213}} =\tau_3.
\end{align*}
Thus, \eqref{st123} holds. We can now conclude from Lemma
\ref{lema:st} that $\sigma_i\trid \tau_i=\tau_i$, $1\leq i \leq
3$, and $\sigma_i\trid \tau_j=\tau_k$, for all $i$, $j$, $k$
distinct. If now $\sigma_3$ acts on \eqref{ts123}, then $
\sigma_3=(\sigma_3 \trid \tau_1) \trid (\sigma_3 \trid
\sigma_2)=\tau_2 \trid \sigma_1$, and \eqref{ts213} holds. \epf

\subsection{Examples of $\D_3^{(2)}$ type}\label{subsubsec:exar2}

We first spell out explicitly Theorem \ref{theorem:dp-cor} and
Corollary \ref{coro:dp-cor} for $p=3$.

\begin{theorem}\label{teor:aplicAHSch}
Let $\sigma_1$, $\sigma_2$, $\sigma_3$, $\tau_1$, $\tau_2$,
$\tau_3\in G$ distinct; denote
$(\sigma,\tau)=(\sigma_1,\sigma_2,\sigma_3,\tau_1,\tau_2,\tau_3)$.
Let $\rho=(\rho,V)\in \widehat{G^{\sigma_1}}$. We assume that
\begin{itemize}
\item[(H1)] $(\sigma,\tau)$ is of type $\D_3^{(2)}$,
\item[(H2)] $(\sigma,\tau)\subseteq \oc$, with $g\in G$ such that $g\trid \sigma_1 =\tau_1$,
\item[(H3)] $q_{\sigma_1\sigma_1}=-1$,
\item[(H4)] there exist $v,w\in V - 0$ such that,
\begin{align}
\label{eqn:d3-a-act} \rho(g^{-1}\sigma_1g) w&=-w,
\\\label{eqn:d3-ga-act}
\rho(\tau_1)v &= -v,
\end{align}
\end{itemize}
Then $\dim\toba(\oc,\rho)=\infty$.\qed
\end{theorem}

\begin{cor}\label{co:especial}
Let $\sigma_1$, $\sigma_2$, $\sigma_3 \in \oc$ distinct. Assume
that there exists $k$, $1\leq k \leq |\sigma_1|$, such that
$\sigma_1^k\neq\sigma_1$ and $\sigma_1^k\in \oc$. Let
$\rho=(\rho,V)\in \widehat{G^{\sigma_1}}$. Assume further that
\begin{itemize}
\item[(1)] $(\sigma_i)_{1\leq i \leq 3}$ is of type $\D_3$,
\item[(2)] $q_{\sigma_1\sigma_1}=-1$.
\end{itemize}
Then $\dim\toba(\oc,\rho)=\infty$.\qed
\end{cor}

Corollary \ref{co:especial} applies notably to a real conjugacy
class of an element of order greater than 2. We list several
applications for $G = \sm$.

\begin{example}\label{exa:12m}
Let $m\geq 6$. Let $\oc$ be the conjugacy class of $\s_m$ of type
$(1^{n_1},2^{n_2},\dots,m^{n_{m}})$, where
\begin{itemize}
    \item $n_1$, $n_2\geq 1$ and
    \item $n_j\geq 1$ for some $j$, $3\leq j \leq m$.
\end{itemize}
Let $\sigma\in \oc$ and $\rho \in \widehat{\s_m^{\sigma}}$. Then
$\dim \toba(\oc,\rho)=\infty$. \end{example}

\pf By hypothesis, we can choose $\sigma = (1\,2) \beta$ where
$\beta$ fixes 1, 2 and 3. If $q_{\sigma\sigma}\neq -1$, then $\dim
\toba(\oc,\rho)=\infty$, by Lemma \ref{odd}. Assume that
$q_{\sigma\sigma}= -1$. Now set
\begin{align*}
x=(1\, 2), \quad y=(1\, 3), \quad z=(2\, 3),\quad \sigma_1=\sigma
= x \beta, \quad\sigma_2=y \beta, \quad \sigma_3:=z \beta.
\end{align*}
Clearly $(\sigma_i)_{1\leq i \leq 3}$ is of type $\D_3$, $\oc$ is
real and $|\sigma_1|>2$. By Corollary \ref{co:especial}, $\dim
\toba(\oc,\rho)=\infty$. \epf

In particular, let $\oc$ be the conjugacy class of $\s_m$ of type
$(1,2,m-3)$, with $m\geq 6$. By the preceding, $\dim
\toba(\oc,\rho)=\infty$. But, if $q_{\sigma\sigma} = -1$, then
$M(\oc, \rho)$ has negative braiding; that is, it is not possible
to decide if the dimension of $\toba(\oc, \rho)$ is infinite via
abelian subracks. See \cite{F-tesis} for details.

\begin{example}\label{co:type2k^3}
Let $m\geq 6$. Let $\sigma\in \sm$ of type
$(1^{n_1},2^{n_2},\dots,m^{n_m})$, $\oc$ the conjugacy class of
$\sigma$ and $\rho \in \widehat{\s_m^{\sigma}}$. Assume that
\begin{itemize}
    \item there
exists $j$, $1\leq j \leq m$, such that $j=2k$, with $k\geq 2$ and
$n_j\geq 3$.
\end{itemize}
Then $\dim \toba(\oc,\rho)=\infty$.
\end{example}

\pf If $q_{\sigma\sigma}\neq -1$, then $\dim
\toba(\oc,\rho)=\infty$, by Lemma \ref{odd}. Assume that
$q_{\sigma\sigma}= -1$. Let
\begin{align*}
\alpha_{1}=(i_{1} \, i_{2}  \, \cdots \, i_{j}), \quad
\alpha_{2}=(i_{j+1} \, i_{j+2} \, \cdots \, i_{2j}), \quad
\alpha_{3}=(i_{2j+1} \, i_{2j+2} \, \cdots \, i_{3j}),
\end{align*}
be three $j$-cycles appearing in the decomposition of $\sigma$ as
product of disjoint cycles and define
\begin{align*}
{\bf I}&=(i_{1}\, i_{3} \, i_{5} \, \cdots \, i_{3j-1}),& B_{1}&=
(i_{1} \,\, i_{j+1})(i_{2} \,\, i_{j+2}) \cdots (i_{j} \,\,
i_{2j}),
\\
{\bf P}&=(i_{2}\, i_{4} \, i_{6} \, \cdots \, i_{3j}),  & B_{2}&=
(i_{j+1} \,\, i_{2j+1})(i_{j+2} \,\, i_{2j+2} ) \cdots ( i_{2j}
\,\, i_{3j} ).\end{align*}

Then  \begin{itemize}
\item[(a)] ${\bf I}$ and ${\bf P}$ are disjoint $3k$-cycles,
\item[(b)] ${\bf I}^k {\bf P}^k=B_{1}B_{2}$,
\item[(c)] $\alpha_{1}\alpha_{2}\alpha_{3}  {\bf I} \alpha_{3}^{-1}\alpha_{2}^{-1}\alpha_{1}^{-1}={\bf P}$,
(and then $\sigma  {\bf I} \sigma^{-1}={\bf P}$),
\item[(d)] ${\bf P}^k \sigma {\bf P}^{k}=\sigma B_1 B_2$, and \item[(e)] ${\bf P}^{-k}\sigma {\bf P}^{-k}=\sigma B_2 B_1$.
\end{itemize}
The first item is clear. To see (b), note that \begin{align*}
B_{1}B_{2}=(i_{1}\,\, i_{j+1}\,\, i_{2j+1})(i_{2} \,\, i_{j+2}\,\,
i_{2j+2} ) \cdots (i_{j} \,\, i_{2j}\,\, i_{3j}).
\end{align*}
(c) follows as in the proof of Lemma \ref{claim:type6k^nj} (c).
(d). By (b) and (c), we have that $\sigma^{-1} {\bf P}^k \sigma
{\bf P}^k={\bf I}^k{\bf P}^k=B_1B_2$, as claimed. (e). By (b) and
(c),  $\sigma^{-1}{\bf P}^{-k}\sigma {\bf P}^{-k}={\bf I}^{-k}
{\bf P}^{-k}=B_2 B_1$ as claimed.

Set now  $\sigma_1:=\sigma$, $\sigma_2:={\bf P}^k \sigma {\bf
P}^{-k}$ and $\sigma_3:={\bf P}^{-k} \sigma {\bf P}^{k}$. As in
the proof of Example \ref{co:type6k^nj} we can see that
$\sigma_1$, $\sigma_2$ and $\sigma_3$ are distinct. We check that
$(\sigma_i)_{1\leq i \leq 3}$ is of type $\D_3$ using Remark
\ref{obs:prop1}.

By (d),  ${\bf P}^k \sigma {\bf P}^k\in \sm^{\sigma}$, i.~e. ${\bf
P}^k \sigma {\bf P}^k \sigma {\bf P}^{-k} \sigma^{-1} {\bf
P}^{-k}=\sigma$, or $\sigma {\bf P}^k \sigma {\bf P}^{-k}
\sigma^{-1}= {\bf P}^{-k}  \sigma {\bf P}^k$. That is,
$\sigma_1\trid \sigma_2= \sigma_3$. Analogously, $\sigma_1\trid
\sigma_3=\sigma_2$ is proved using (e). To check that $\sigma_2
\trid \sigma_3 =\sigma_1$, note that $\sigma_2\trid \sigma_3={\bf
P}^k \sigma {\bf P}^{-k}  {\bf P}^{-k} \sigma {\bf P}^{k}  {\bf
P}^k \sigma^{-1} {\bf P}^{- k}=\sigma$, because ${\bf P}^k \sigma
{\bf P}^{-2k}={\bf P}^k \sigma {\bf P}^{k}{\bf P}^{-3k}=\sigma
B_1B_2 \in \sm^{\sigma}$, by (a) and (d).

We now apply Corollary \ref{co:especial} and conclude that $\dim
\toba(\oc,\rho)=\infty$. \epf

\medbreak We shall need a few well-known results on symmetric
groups.

\begin{obs}\label{obs:repfielSn}
(i) If $\rho$ is a faithful representation of $\s_n$, then
$\rho(\tau)\notin \ku\Id$, for every $\tau \in \s_n$, $\tau\neq
\id$ (since $\s_n$ is centerless).

(ii) If $\rho=(\rho,W)\in \widehat{\s_n}$, with $\rho\neq \sgn$,
then for any involution $\tau\in\s_n$ (i.~e., $\tau^2=\id$), there
exists $w\in W-0$ such that $\rho(\tau)w=w$ (otherwise $\rho(\tau)
=-\Id$).
\end{obs}

\begin{example}\label{co:type2^3}
Let $m\geq 6$. Let $\sigma\in \sm$ of type
$(1^{n_1},2^{n_2},\dots,m^{n_m})$, $\oc$ the conjugacy class of
$\sigma$ and $\rho \in \widehat{\s_m^{\sigma}}$. Assume that
\begin{itemize}
    \item $n_2\geq 3$ and
    \item there exists $j$, with $j\geq 3$, such that
$n_j\geq 1$.
\end{itemize}

Then $\dim \toba(\oc,\rho)=\infty$.
\end{example}

\pf By Lemma \ref{odd}, we may suppose that $q_{\sigma\sigma}=
-1$. Assume that $(i_1\, i_2)$, $(i_3\, i_4)$ and $(i_5\, i_6)$
are three transpositions appearing in the decomposition of
$\sigma$ as a product of disjoint cycles. We define
\begin{align*}
x:=(i_1\, i_2)(i_3\, i_4)(i_5\, i_6),\quad  y:=(i_1\, i_4)(i_3\,
i_6)(i_2\, i_5), \quad z:=(i_1\, i_6)(i_2\, i_3)(i_4\, i_5)
\end{align*}
and $\alpha:=x \sigma$. It is easy to see, using for instance
Proposition \ref{prop:caracR2}, that
\begin{align*}
\sigma_1:=\sigma, \qquad \sigma_2:=y \alpha, \qquad \sigma_3:=z
\alpha,
\end{align*}
is of type $\D_3$. Then $\dim \toba(\oc,\rho)=\infty$, by
Corollary \ref{co:especial}. Indeed, $\sigma^{-1}\in \Oc$, but
$\sigma\neq \sigma^{-1}$ because $\sigma$ has order $>2$. \epf

In the proof of the next Example, we need some notation for the
induced representation. Let $H$ be a subgroup of a finite group
$G$ of index $k$, $\phi_1, \dots, \phi_k$ the left cosets of $H$
in $G$, with representatives $g_{\phi_1}, \dots, g_{\phi_k}$. Let
$\theta=(\theta,W)\in \widehat{H}$, and $w_1, \dots w_r$ a basis
of $W$. Set $V:=$span-$\{g_{\phi_i} w_j \, | \, 1\leq i \leq k,
1\leq j \leq r\}$. For $i$, $j$, with $1\leq i \leq k$, $1\leq j
\leq r$ we define $\rho:G \to \Aut(V)$ by
\begin{align}\label{eq:repind}
\rho(g)(g_{\phi_i}w_j)= g_{\phi_{l}} \theta(h)w_j, \qquad \text{
where $gg_{\phi_i}=g_{\phi_{l}}h$, with $h\in H$.}
\end{align}
Thus $\rho=(\rho,V)$ is a representation of $G$ and
$\deg\rho=[G:H]\deg\theta$.

\begin{example}\label{co:type2n2:6}
Let $m\geq 12$. Let $\sigma\in \sm$ of type
$(1^{n_1},2^{n_2},\dots,m^{n_m})$, $\oc$ the conjugacy class of
$\sigma$ and $\rho \in \widehat{\s_m^{\sigma}}$. If $n_2\geq 6$,
then $\dim \toba(\oc,\rho)=\infty$.
\end{example}

\pf By Lemma \ref{odd}, we may suppose that $q_{\sigma\sigma}=
-1$. We denote the $n_2$ transpositions appearing in the
decomposition of $\sigma$ as product of disjoint cycles by
$A_{1,2},\dots,A_{n_2,2}$ and we define $A_2=A_{1,2}\cdots
A_{n_2,2}$. Let us suppose that $A_{1,2}=(i_1\, i_2)$,
$A_{2,2}=(i_3\, i_4)$, $A_{3,2}=(i_5\, i_6)$, $A_{4,2}=(i_7\,
i_8)$, $A_{5,2}=(i_9\, i_{10})$ and $A_{6,2}=(i_{11}\, i_{12})$.
We define $x:=(i_1\, i_2)(i_3\, i_4)(i_5\, i_6)(i_7\, i_8)(i_9\,
i_{10})(i_{11}\, i_{12})$ and $\alpha:=x \sigma$.

If there exists $j$, with $j\geq 3$, such that $n_j\geq 1$, then
the result follows from Example \ref{co:type2^3}. Assume that
$n_j=0$, for every $j\geq 3$, i.~e. the type of $\sigma$ is
$(1^{n_1}, 2^{n_2})$. The centralizer of $\sigma$ in $\s_m$ is
$\s_m^{\sigma}=T_1\times T_2$, with $T_1\simeq \s_{n_1}$ and
$T_2=\Gamma\rtimes \Lambda$, with
\begin{align*}
\Gamma:=\langle A_{1,2},\dots,A_{n_2,2} \rangle, \quad
\Lambda:=\langle B_{1,2},\dots,B_{n_2-1,2} \rangle.
\end{align*}
Here $B_{l,2} := (i_{2l-1}\,\, i_{2l+1})(i_{2l}\,\, i_{2l+2})$, for $1\le l \le n_2-1$.
Note that $\Gamma \simeq(\Z/2)^{n_2}$ and $\Lambda\simeq\s_{n_2}$.
Now, $\rho=\rho_1 \otimes \rho_2$, with $\rho_1=(\rho_1,V_1) \in
\widehat{T_1}$ and $\rho_2=(\rho_2,V_2) \in \widehat{T_2}$.

For every $t$, $1\leq t \leq n_2$, we define $\chi_t\in
\widehat{\Gamma}$, by $\chi_{t}(A_{l,2})=(-1)^{ \delta_{t,l}}$,
$1\leq l \leq n_2$. Then, the irreducible representations of
$\Gamma$ are
$$\chi_{t_1,...,t_J}:= \chi_{t_1}\dots\chi_{t_J},
\qquad 0 \leq J \leq n_2,\quad 1 \leq t_1 < \cdots <t_J \leq
n_2.$$ The case $J=0$ corresponds to the trivial representation of
$\Gamma$.

For every $J$, with $0 \leq J \leq n_2$, we denote
$\chi_{(J)}:=\chi_{1,\dots,J}$. The action of $\Lambda$ on
$\Gamma$ induces a natural action of $\Lambda$ on
$\widehat{\Gamma}$, namely $(\lambda \cdot
\chi)(A_{l,2}):=\chi(\lambda^{-1}A_{l,2}\lambda)$, $1\leq l \leq
n_2$, $\lambda\in \Lambda$. The orbit and the isotropy subgroup of
$\chi_{(J)} \in \widehat{\Gamma}$ are
\begin{align}
\Oc_{\chi_{(J)}}\label{eq:orb}&=\{\chi_{k_1,\dots,k_J} : 1 \leq
k_1 <\dots<k_J \leq n_2\},\\
\Lambda^{\chi_{(J)}}\label{eq:iso}&=(\Lambda^{\chi_{(J)}})_1\times
(\Lambda^{\chi_{(J)}})_2 \\ \notag &=\langle
B_{1,2},\dots,B_{J-1,2} \rangle \times \langle
B_{J+1,2},\dots,B_{n_2-1,2} \rangle \simeq \s_{J} \times
\s_{n_2-J}.
\end{align}
Thus, the characters $\chi_{(J)}$, $0 \leq J \leq n$, form a
complete set of representatives of the orbits in
$\widehat{\Gamma}$ under the action of $\Lambda$.

Since $\rho_2 \in \widehat{\Gamma \rtimes\Lambda}$, we have that
$\rho_2=\Ind_{\Gamma \rtimes\Lambda^{\chi_{(J)}}}^{\Gamma
\rtimes\Lambda}\chi_{(J)}\otimes \mu$, with $\chi_{(J)}$ as above
and $\mu=(\mu,W) \in \widehat{\Lambda^{\chi_{(J)}}}$ -- see
\cite[Section 8.2]{S}. By \eqref{eq:iso}, $\mu=\mu_1\otimes
\mu_2$, with $\mu_l=(\mu_l,W_l) \in
\widehat{(\Lambda^{\chi_{(J)}})_l}$, $l=1$, $2$. Let
$\{\phi_1=\Lambda^{\chi_{(J)}},\dots,\phi_k\}$ the left cosets of
$\Lambda^{\chi_{(J)}}$ in $\Lambda$, where
$k=[\Lambda:\Lambda^{\chi_{(J)}}]=\frac{n_2!}{J!(n_2-J)!}$.

Note that $$B_{1,2}= (i_1\,\, i_3)(i_2\,\, i_4), \quad
B_{3,2}=(i_5\,\, i_7)(i_6\,\, i_8)\quad \text{ and }\quad
B_{5,2}=(i_9\,\, i_{11})(i_{10}\,\, i_{12}).$$ We define
$B:=B_{1,2}B_{3,2}B_{5,2}$. Notice that the order of $B$ is $2$.

Since $q_{\sigma\sigma}=-1$, then $J$ is odd. We will consider two
cases.

\emph{{\bf CASE (1):} assume that $J\leq 5$}. Then, $B\not\in
\Lambda^{\chi_{(J)}}$. This implies that the left coset  $\phi$ of
$\Lambda^{\chi_{(J)}}$ in $\Lambda$ containing $B$ is not the
trivial coset $\phi_1$. We choose as representatives of the cosets
$\phi_1$ and $\phi$ to $g_{\phi_1}=\id$ and $g_{\phi}=B$,
respectively. We define $v_2:=g_{\phi_1}w+g_{\phi}w$, with $w\in
W-0$. Notice that $B g_{\phi_1}=g_{\phi} \id$ and $B
g_{\phi}=g_{\phi_1} \id$. Using \eqref{eq:repind}, we have that
\begin{equation}\label{eq:rho2v2}
\begin{aligned}
\rho_2(B)v_2&=\rho_2(B)(g_{\phi_1}w)+\rho_2(B)(g_{\phi}w)\\&=g_{\phi}
\mu(\id)w+g_{\phi_1} \mu(\id)w=g_{\phi}w+g_{\phi_1}w=v_2.
\end{aligned}
\end{equation}

Let $v:=v_1\otimes v_2$, with $v_1\in V_1-0$. Then,
\begin{align}\label{eq:rho2v}
\rho(B)v=(\rho_1\otimes\rho_2)(\id,B)(v_1\otimes
v_2)=\rho_1(\id)v_1\otimes \rho_2(B)v_2=v_1\otimes v_2=v,
\end{align}
by \eqref{eq:rho2v2}. We define $\sigma_1:=\sigma$,
\begin{align*}
\sigma_2&:=(i_1\,\, i_6)(i_3\,\, i_8)(i_5\,\, i_{10})(i_7\,\, i_{12})(i_{9}\,\, i_2)(i_{11}\,\, i_4)\alpha, \\
\sigma_3&:=(i_1\,\, i_{10})(i_3\,\, i_{12})(i_5\,\, i_2)(i_7\,\, i_4)(i_9\,\, i_6)(i_{11}\,\, i_8)\alpha,\\
\tau_1&:=(i_1\,\, i_4)(i_3\,\, i_2)(i_5\,\, i_8)(i_7\,\, i_6)(i_9\,\, i_{12})(i_{11}\,\, i_{10})\alpha, \\
\tau_2&:=(i_1\,\, i_8)(i_3\,\, i_6)(i_5\,\, i_{12})(i_7\,\, i_{10})(i_9\,\, i_4)(i_{11}\,\, i_2)\alpha,\\
\tau_3&:=(i_1\,\, i_{12})(i_3\,\, i_{10})(i_5\,\, i_4)(i_7\,\,
i_2)(i_9\,\, i_8)(i_{11}\,\, i_6)\alpha.
\end{align*}
We can check by straightforward computations that $(\sigma,\tau)$
is of type $\D_3^{(2)}$. 
Let $g:=(i_2\,\, i_4)(i_6\,\, i_8)(i_{10}\,\, i_{12})$; thus,
$g\trid \sigma=\tau_1$. Moreover, $\tau_1=\sigma B=g\sigma g$ and
$\sigma_2\tau_2=B=g\sigma_2\tau_2 g$. Then,
\begin{align*}
\rho(\tau_1)v=-v=\rho(g\sigma_1 g)v,
\end{align*}
by \eqref{eq:rho2v}. Therefore, $\dim \toba(\oc,\rho)=\infty$, by
Theorem \ref{teor:aplicAHSch}.

\smallbreak

\emph{{\bf CASE (2):} assume that $J\geq 7$}. Then, $B\in
\Lambda^{\chi_{(J)}}$; moreover, $B\in (\Lambda^{\chi_{(J)}})_1$.
Also, $B g_{\phi_1}=g_{\phi_1} B$.

Let $v_2=g_{\phi_1} w$, with $w\in W-0$. Since $W=W_1\otimes W_2$,
we may assume that $w=w_1\otimes w_2$, with $w_1\in W_1-0$ and $w_2\in
W_2-0$. Then, using \eqref{eq:repind},
\begin{align*}
\rho_2(B)v_2&=\rho_2(B)(g_{\phi_1}w)=g_{\phi_1}
\mu(B)w=g_{\phi_1}(\mu_1\otimes\mu_2)(B,\id)(w_1\otimes w_2)\\
&=g_{\phi_1}\Big(\mu_1(B)(w_1)\otimes
\mu_2(\id)(w_2)\Big)=g_{\phi_1}\Big((\mu_1(B)(w_1)\otimes
w_2\Big).
\end{align*}
Notice that $\mu_1\in \widehat{(\Lambda^{\chi_{(J)}})_1}$. Since
$(\Lambda^{\chi_{(J)}})_1\simeq \s_{J}$, if $\mu_1\neq \sgn$, with
$\sgn$ the sign representation of $\s_{J}$, then there exists
$w_1\in W_1-0$ such that $\mu_1(B)(w_1)=w_1$, by Remark
\ref{obs:repfielSn} (ii). In this case, we have
\begin{align}\label{eq:rho2v2Caso2}
\rho_2(B)v_2=g_{\phi_1}(\mu_1(B)(w_1)\otimes
w_2)=g_{\phi_1}(w_1\otimes w_2)=g_{\phi_1}w=v_2.
\end{align}
Taking $v:=v_1\otimes v_2$, with $v_1\in V_1-0$, we have
\begin{align*}
\rho(B)v=(\rho_1\otimes\rho_2)(\id,B)(v_1\otimes
v_2)=\rho_1(\id)v_1\otimes \rho_2(B)v_2=v_1\otimes v_2=v,
\end{align*}
by \eqref{eq:rho2v2Caso2}. Considering $\sigma_i$, $\tau_i$,
$1\leq i \leq 3$, as in the previous case, the hypothesis of
Corollary \ref{co:especial} hold. Therefore, $\dim
\toba(\oc,\rho)=\infty$.

On the other hand, let us suppose that $\mu_1= \sgn$. Let $w\in
W$, with $w=w_1\otimes w_2$, $w_1\in W_1-0$ and $w_2\in W_2-0$.
Let $v_2=g_{\phi_1} w$; since $\mu_1(B)(w_1)=-w_1$, we have
$\rho_2(B)v_2=-v_2$. Choosing $v:=v_1\otimes v_2$, with $v_1\in
V_1-0$, we have that
\begin{align}\label{eq:caso:sign}
\rho(B)v=(\rho_1\otimes\rho_2)(\id,B)(v_1\otimes
v_2)=\rho_1(\id)v_1\otimes \rho_2(B)v_2=- v_1\otimes v_2=-v.
\end{align}
We define $\overline{\sigma_1}:=\sigma$,
\begin{align*}
\overline{\sigma_2}&:=(i_1\,\, i_6)(i_4\,\, i_7)(i_5\,\, i_{10})(i_8\,\, i_{11})(i_{2}\,\, i_9)(i_{3}\,\, i_{12})\alpha, \\
\overline{\sigma_3}&:=(i_1\,\, i_{10})(i_4\,\, i_{11})(i_2\,\, i_5)(i_3\,\, i_8)(i_6\,\, i_9)(i_{7}\,\, i_{12})\alpha,\\
\overline{\tau_1}&:=(i_1\,\, i_3)(i_2\,\, i_4)(i_5\,\, i_7)(i_6\,\, i_8)(i_9\,\, i_{11})(i_{10}\,\, i_{12})\alpha, \\
\overline{\tau_2}&:=(i_1\,\, i_7)(i_2\,\, i_{12})(i_3\,\, i_{9})(i_4\,\, i_{6})(i_5\,\, i_{11})(i_{8}\,\, i_{10})\alpha,\\
\overline{\tau_3}&:=(i_1\,\, i_{11})(i_2\,\, i_{8})(i_3\,\,
i_5)(i_4\,\, i_{10})(i_6\,\, i_{12})(i_{7}\,\, i_9)\alpha.
\end{align*}
It can be shown that $(\overline{\sigma},\overline{\tau})$ is of
type $\D_3^{(2)}$.
Let now $\overline{g}= (i_2\,\, i_3)(i_6\,\,
i_7)(i_{10}\,\, i_{11})$; then, $\overline{g}\trid
\sigma=\overline{\tau_1}$. Furthermore, $\overline{\tau_1}=
B=\overline{g}\sigma \overline{g}$ and $\overline{\sigma_2} \,
\overline{\tau_2}=\sigma B=\overline{g}\, \overline{\sigma_2}\,
\overline{\tau_2} \, \overline{g}$. Then
\begin{align*}
\rho(\overline{\tau_1})v=-v=\rho(\overline{g}\sigma \overline{g})v
\qquad \text{ and } \qquad \rho(\overline{\sigma_2}\,
\overline{\tau_2})v=v=\rho(\overline{g}\, \overline{\sigma_2}\,
\overline{\tau_2}\, \overline{g})v,
\end{align*}
by \eqref{eq:caso:sign}. Therefore, $\dim \toba(\oc,\rho)=\infty$,
by Theorem \ref{teor:aplicAHSch}.\epf

\medbreak

A way to obtain a family of type $\D_3$ is to start from a
monomorphism $\rho:\mathbb S_3\to G$ and to consider the image
by $\rho$ of the transpositions. Another way is as follows.

\begin{rmk}\label{exa:triple-central}Let $G$ be a finite group and $z\in Z(G)$.

(i). Let $(\sigma_i)_{i\in \Z/3}$ be of type $\D_3$. Then
$(z\sigma_i)_{i\in \Z/3}$ is also of type $\D_3$.

(ii). Let $(\sigma, \tau) = (\sigma_i)_{i\in \Z/3}\cup
(\tau_i)_{i\in \Z/3}$ be a family of type $\D_3^{(2)}$. Then
$(z\sigma, z\tau) = (z\sigma_i)_{i\in \Z/3}\cup (z\tau_i)_{i\in
\Z/3}$ is also a family of type $\D_3^{(2)}$.
\end{rmk}

Here is a combination of these two ways.

\begin{exa}\label{exa:gl2-triple} Let $p$ be a prime number and $q =
p^m$, $m\in \N$, such that 3 divides $q-1$. Let $\omega \in \kc$
be a primitive third root of 1.

(i). If $c\in \kc$, then $(\mu_i)_{i\in \Z/3}$, where
$\mu_i = \begin{pmatrix}0 & \omega^i \\ \omega^{2i}c & 0
\end{pmatrix}$, is a family of type $\D_3$ in $\mathbf{GL}(2, \kc)$. If $c=-1$, then this
is a family of type $\D_3$ in $\mathbf{SL}(2, \kc)$. The orbit of
$\mu_i$ is the set of matrices with minimal polynomial $T^2 -c$.

(ii). Let $N > 3$ be an integer and let $\Tb$ be the subgroup of
diagonal matrices in $\mathbf{GL}(N, \kc)$. Let $ \lambda =
\diag(\lambda_1, \lambda_2,\dots , \lambda_N) \in \Tb$. Let $\Oc$
be the conjugacy class of $\lambda$. Assume that $\lambda_1 =
-\lambda_2$ and let $c = \lambda_1^2$. Assume also that there
exist $i,j$, with $3\le i,j\le N$ such that $\lambda_i \neq
\lambda_j$; say $i=3$, $j=4$, for simplicity of the exposition.
Then  $(\sigma_i)_{i\in \Z/3} \cup (\tau_i)_{i\in \Z/3}$, where
$$
\sigma_i= \begin{pmatrix} \mu_i& 0  \\
0 &\diag(\lambda_3, \lambda_4,\dots , \lambda_N)
\end{pmatrix},\qquad \tau_i= \begin{pmatrix} \mu_i& 0  \\
0 &\diag(\lambda_4, \lambda_3,\dots , \lambda_N)
\end{pmatrix},
$$
is a family of type $\D_3^{(2)}$ in the orbit $\Oc\subset
\mathbf{GL}(N, \kc)$. \end{exa}Let $\W = \s_N$ act on $\Tb$ in the
natural way.  Let $\chi:\mathbf{GL}(N, \kc)\to \ku^{\times}$ be a
character; it restricts to an irreducible representation $(\chi,
\ku)$ of the centralizer $\mathbf{GL}(N, \kc)^{\sigma_0}$. Fix a
group isomorphism $\varphi: \kc^{\times}\to \G_{q-1} \subset
\ku^{\times}$, where $\G_{q-1}$ is the group of $(q-1)$-th roots
of 1 in $\ku$. Recall that $\chi = \varphi(\det^h)$ for some
integer $h$. Thus the restriction of $\chi$ to $\Tb$ is
$\W$-invariant.

\begin{prop}\label{prop:r2-gln} Keep the notation above.  Assume that
$\chi(\lambda) = -1$. Then the dimension of the Nichols algebra
$\toba(\Oc, \chi)$ is infinite.
\end{prop}

\pf The result follows from Theorem \ref{teor:aplicAHSch}. Indeed,
hypothesis (H1) and (H2) clearly hold. The matrix $g =
\begin{pmatrix} \id_2& 0&0 &0
\\ 0 &0 &1&0
\\ 0 &1 &0&0
\\0 &0 &0& \id_{N-4}
\end{pmatrix}$ is an involution that satisfies $g\trid \sigma_0 =
\tau_0$. Because of the explicit form of $\chi$, $\chi(\sigma_0) =
-1 = \chi(\tau_0)$, hence (H3) and (H4) hold. \epf

This example can be adapted to the setting of semisimple orbits in
finite groups of Lie type.


\section{A technique from the symmetric group $\sk$}

The classification of the finite-dimensional Nichols algebras over
$\sk$, given in \cite{AHS}, relies on the fact (proved in
\emph{loc.~cit.}) that some Nichols algebras $\toba(V_i\oplus
V_j)$ have infinite dimension. According to the general strategy
proposed in the present paper, each of these pairs $(V_i,V_j)$
gives rise to a rack and a cocycle, and to a technique to discard
Nichols algebras over other groups. Here we study one of these
possibilities, and leave the others for a future publication.

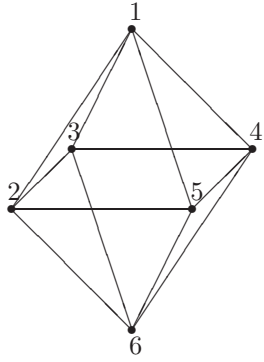
\begin{figure}[ht]
\vspace{2cm}
\begin{align*}
\setlength{\unitlength}{0.8cm}
\begin{picture}(3.5,0)
\put(0,0){\circle*{.15}} \put(1,1){\circle*{.15}}
\put(3,0){\circle*{.15}} \put(4,1){\circle*{.15}}
\put(2,3){\circle*{.15}} \put(2,-2){\circle*{.15}}
\put(0,0){\line(2,3){2}} \put(0,0){\line(1,-1){2}}
\put(0,0){\line(1,1){1}} \put(0,0){\line(3,0){3}}
\put(1,1){\line(3,0){3}} \put(1,1){\line(1,2){1}}
\put(1,1){\line(1,-3){1}} \put(2,3){\line(1,-1){2}}
\put(3,0){\line(1,1){1}} \put(3,0){\line(-1,3){1}}
\put(3,0){\line(-1,-2){1}} \put(2,-2){\line(2,3){2}}
\put(-0.07,0.15){\small{2}} \put(2.97,0.15){\small{5}}
\put(0.93,1.15){\small{3}} \put(3.95,1.15){\small{4}}
\put(1.95,3.15){\small{1}} \put(1.95,-2.4){\small{6}}
\end{picture}
\end{align*}
\vspace{1.4cm}\caption{Octahedral rack.}\label{rackocta}
\end{figure}

\medbreak The \emph{octahedral rack} is the rack
$X=\{1,2,3,4,5,6\}$ given by the vertices of the octahedron with
the operation of rack given by the ``right-hand rule'', i.~e. if
$T_i$ is the orthogonal linear map that fixes $i$ and rotates the
orthogonal plane by an angle of $\pi/2$ with the right-hand rule
(pointing the thumb to $i$), then we define $\trid:X\times X \to
X$ by $i\trid j:=T_i(j)$ -- see Figure \ref{rackocta}.

Explicitly,
\begin{align*}
1\trid 1=1, \quad 2\trid 1=3, \quad 3\trid 1=4, \quad 4\trid 1=5,
\quad 5\trid 1=2, \quad 6\trid 1=1,\\
1\trid 2=5, \quad 2\trid 2=2, \quad 3\trid 2=1, \quad 4\trid 2=2,
\quad 5\trid 2=6, \quad 6\trid 2=3,\\
1\trid 3=2, \quad 2\trid 3=6, \quad 3\trid 3=3, \quad 4\trid 3=1,
\quad 5\trid 3=3, \quad 6\trid 3=4,\\
1\trid 4=3, \quad 2\trid 4=4, \quad 3\trid 4=6, \quad 4\trid 4=4,
\quad 5\trid 4=1, \quad 6\trid 4=5,\\
1\trid 5=4, \quad 2\trid 5=1, \quad 3\trid 5=5, \quad 4\trid 5=6,
\quad 5\trid 5=5, \quad 6\trid 5=2,\\
1\trid 6=6, \quad 2\trid 6=5, \quad 3\trid 6=2, \quad 4\trid 6=3,
\quad 5\trid 6=4, \quad 6\trid 6=6.
\end{align*}

\smallbreak

Let $G$ be a finite group, $\sigma_1$, $\sigma_2$, $\sigma_3$,
$\sigma_4$, $\sigma_5$, $\sigma_6\in G$ distinct elements and
$\oc$ the conjugacy class of $\sigma_1$ in $G$.

\begin{definition}\label{de:oct}
We will say that $(\sigma_i)_{1\leq i \leq 6}$ is of \emph{type
$\oct$} if the following holds
\begin{align*}
\sigma_i \trid \sigma_j=\sigma_{i\trid j}, \qquad 1\leq i,j \leq
6.
\end{align*}
\emph{Here and in the rest of this section, $\trid$ in the
subindex is the operation of rack in the octahedral rack.} In
other words, $(\sigma_i)_{1\leq i \leq 6}$ is of type $\oct$ if
and only if $\{\sigma_i \, | \, 1\leq i \leq 6 \}$ is isomorphic
to the octahedral rack via $i\mapsto \sigma_i$.
\end{definition}

\begin{example}
Let $m\geq 4$. Let us consider in $\s_m$ the following 4-cycles
\begin{equation}\label{eq:relss4}\begin{aligned}
\widetilde{\sigma}_1&=(1\,2\,3\,4), \quad
&\widetilde{\sigma}_2&=(1\,2\,4\,3), &\quad
\widetilde{\sigma}_3&=(1\,3\,2\,4),\\
\widetilde{\sigma}_4&=(1\,3\,4\,2), \quad
&\widetilde{\sigma}_5&=(1\,4\,2\,3), \quad
&\widetilde{\sigma}_6&=(1\,4\,3\,2).
\end{aligned}\end{equation}
It is easy to see that $(\widetilde{\sigma}_i)_{1\leq i \leq 6}$
satisfy the relations given in the previous definition. Thus,
$(\widetilde{\sigma}_i)_{1\leq i \leq 6}$ is of  type $\oct$.
\end{example}

\smallbreak Let $\chi_-\in \widehat{\s_4^{\widetilde{\sigma}_1}}$
be given by $\chi_-(1\,2\,3\,4)=-1$. The goal of this Section is
to apply the next result, cf. \cite[Theor. 4.7]{AHS}.

\begin{theorem}\label{teor:AHS2}
The Nichols algebra $\toba\left(M(\oc_4^4, \chi_-)\oplus
M(\oc_4^4, \chi_-)\right)$ has infinite dimension. \qed
\end{theorem}

\begin{remark}\label{obs:isoEVTocta}
We note that $M(\oc_4^4, \chi_-)\oplus M(\oc_4^4, \chi_-)\simeq
(\ku Y, \q)$ as braided vector spaces, where
\begin{itemize}
\item $Y=\{x_i, y_j\ | \, 1\leq i,j\leq 6\}\simeq X^{(2)}$, see
Definition \ref{def:square-rack};
\item $\q$ is the constant cocycle $\q \equiv -1$.
\end{itemize}
\pf We define
\begin{align*}
\widetilde{\sigma}_1&:=(1\,2\,3\,4)=:\widetilde{\tau}_1, \quad
&\widetilde{\sigma}_2&:=(1\,2\,4\,3)=:\widetilde{\tau}_2, \quad
&\widetilde{\sigma}_3&:=(1\,3\,2\,4)=:\widetilde{\tau}_3,\\
\widetilde{\sigma}_4&:=(1\,3\,4\,2)=:\widetilde{\tau}_4, \quad
&\widetilde{\sigma}_5&:=(1\,4\,2\,3)=:\widetilde{\tau}_5, \quad
&\widetilde{\sigma}_6&:=(1\,4\,3\,2)=:\widetilde{\tau}_6.
\end{align*}
We will denote by $(\widetilde{\sigma}_j)_{1 \leq j \leq 6}$
(resp. $(\widetilde{\tau}_j)_{1 \leq j \leq 6}$) the first copy
(resp. the second copy) of $\oc_4^4$, with system of left cosets
representatives of $\s_4^{(1\,2\,3\,4)}$ given by
$\widetilde{g}_1=\widetilde{g}_{7}=\widetilde{\sigma}_1$,
$\widetilde{g}_2=\widetilde{g}_{8}=\widetilde{\sigma}_5$,
$\widetilde{g}_3=\widetilde{g}_{9}=\widetilde{\sigma}_2$,
$\widetilde{g}_4=\widetilde{g}_{10}=\widetilde{\sigma}_3$,
$\widetilde{g}_5=\widetilde{g}_{11}=\widetilde{\sigma}_4$,
$\widetilde{g}_6=\widetilde{g}_{12}=\widetilde{\sigma}_2^2\widetilde{\sigma}_1$.
The map $M(\oc_4^4, \chi_{-})\oplus M(\oc_4^4, \chi_{-}) \to (\ku Y, \mathfrak q)$ given by
\begin{align*}
\widetilde{g}_i \mapsto x_i \quad \text{ and } \quad
\widetilde{g}_{i+6} \mapsto y_i, \qquad 1\leq i \leq 6,
\end{align*}
is an isomorphism of braided vector spaces. \epf
\end{remark}

\begin{prop}
A family $(\sigma_i)_{1\leq i \leq 6}$ of  distinct elements in
$G$ is of type $\oct$ if and only if the following identities
hold:
\begin{align}
\label{rel1:R4*}\sigma_1\trid \sigma_2&=\sigma_5,  &\sigma_1\trid
\sigma_3&=\sigma_2,  &\sigma_1\trid \sigma_4&=\sigma_3,
&\sigma_1\trid \sigma_5&=\sigma_4, &\sigma_1\trid \sigma_6&=\sigma_6,\\
\label{rel2:R4*}\sigma_2\trid \sigma_1&=\sigma_3, &\sigma_2\trid
\sigma_3&=\sigma_6,  &\sigma_2\trid \sigma_4&=\sigma_4,
&\sigma_2\trid \sigma_5&=\sigma_1, &\sigma_2\trid
\sigma_6&=\sigma_5.
\end{align}
\end{prop}
\pf If we apply $\sigma_1\trid\underline{\quad}$ to the relations
in \eqref{rel2:R4*}, then we obtain the relations
$\sigma_5\trid\sigma_j=\sigma_{5\trid j}$, $1\leq j \leq 6$,
because $\sigma_1\trid \sigma_2=\sigma_5$. Analogously, we obtain
the relations $\sigma_i\trid\sigma_j=\sigma_{i\trid j}$, $1\leq j
\leq 6$, for $i=3$, $4$; and the relations
$\sigma_6\trid\sigma_j=\sigma_{6\trid j}$, $1\leq j \leq 6$,
follow by applying $\sigma_5\trid \underline{\quad}$ to the ones
in \eqref{rel2:R4*}. \epf

\begin{lema}\label{le:^4iguales}
If $(\sigma_i)_{1\leq i \leq 6}$ is of type $\oct$,
then
\begin{itemize}
\item[(i)]
$\sigma_1^4=\sigma_2^4=\sigma_3^4=\sigma_4^4=\sigma_5^4=\sigma_6^4$,
\item[(ii)] $\sigma_1\sigma_6=\sigma_2\sigma_4=\sigma_3\sigma_5$,
\item[(iii)]
$\sigma_2^2\sigma_5^2=\sigma_1^3\sigma_6=\sigma_3^2\sigma_2^2$,
\item[(iv)]
$\sigma_5^2\sigma_2^2=\sigma_1\sigma_6^3=\sigma_2^2\sigma_3^2$.
\end{itemize}
\end{lema}

\pf (i). Since
$\sigma_i\trid(\sigma_i\trid(\sigma_i(\trid(\sigma_i\trid
\sigma_j))))=\sigma_j$, then $\sigma_i^4\in G^{\sigma_j}$, $1\leq
i,j \leq 6$. Hence $\sigma_1^4=
(\sigma_3\sigma_2\sigma_3^{-1})^4=\sigma_3\sigma_2^4\sigma_3^{-1}=
\sigma_2^4$, and the rest is similar. (ii). By Definition
\ref{de:oct}, we see that \begin{align*} \sigma_3\sigma_5=\sigma_3
\sigma_1\sigma_2\sigma_1^{-1}= \sigma_3
\sigma_2\sigma_5\sigma_2^{-1} \sigma_2\sigma_1^{-1}=
\sigma_2 \sigma_1\sigma_5\sigma_1^{-1}=\sigma_2 \sigma_4,\\
\sigma_3\sigma_5=\sigma_3   \sigma_2\sigma_6\sigma_2^{-1}=
\sigma_3 \sigma_6\sigma_5\sigma_6^{-1} \sigma_6\sigma_2^{-1}=
\sigma_6 \sigma_2\sigma_5\sigma_2^{-1}=\sigma_6 \sigma_1.
\end{align*}
Then, $\sigma_1\sigma_6=\sigma_2\sigma_4=\sigma_3\sigma_5$, as
claimed.

(iii). By (ii), we have that
\begin{align*}
\sigma_2^2\sigma_5^2=\sigma_2 \sigma_5 \sigma_1  \sigma_5=\sigma_5
\sigma_1 \sigma_1  \sigma_5=\sigma_5 \sigma_1  \sigma_4 \sigma_1=
\sigma_5  \sigma_3 \sigma_1^2= \sigma_1  \sigma_6
\sigma_1^2=\sigma_1^3  \sigma_6.
\end{align*}
Then, $\sigma_2^2\sigma_5^2=\sigma_1^3  \sigma_6$. We apply
$\sigma_1 \trid (\sigma_1 \trid(\sigma_1 \trid \underline{\quad}
))$ to the last expression and we have
$\sigma_3^2\sigma_2^2=\sigma_1^3 \sigma_6$.

(iv) follows from (iii) applying $\sigma_2 \trid(\sigma_2
\trid\underline{\quad})$. \epf

\begin{definition}
Let $\sigma_i$, $\tau_i\in G$, $1\leq i \leq 6$, all distinct. We
say  that $(\sigma,\tau)$ is \emph{of type $\oct^{(2)}$} if
$(\sigma_i)_{1\leq i \leq 6}$ and $(\tau_j)_{1\leq j \leq 6}$ are
both of type $\oct$, and
\begin{align}\label{R4^2}
\sigma_i\trid \tau_j=\tau_{i\trid j}, \qquad  \tau_i\trid
\sigma_j=\sigma_{i\trid j}, \qquad 1\leq i,j \leq 6.
\end{align}
\end{definition}

\begin{lema}\label{le:sigma:tau}
If $(\sigma,\tau)$ is of type $\oct^{(2)}$, then
\begin{itemize}
\item[(i)]
$\sigma_1\tau_6=\sigma_6\tau_1=\sigma_2\tau_4=\sigma_4\tau_2=\sigma_3\tau_5=\sigma_5\tau_3$,
\item[(ii)] $\sigma_j^{-1}\tau_j=\sigma_1^{-1}\tau_1$, $2 \leq j \leq
6$, 
\item[(iii)] $\tau_2^{-2}\sigma_5\tau_5=\tau_1^{-1}\sigma_6$,
\item[(iv)] $\tau_2^{-2}\sigma_3\tau_3=\sigma_1\tau_6^{-1}$,
\item[(v)] $\sigma_2^{-2}\sigma_5\tau_5=\sigma_1^{-2}\tau_1\sigma_6$,
\item[(vi)] $\sigma_2^{-2}\sigma_3\tau_3=\tau_1\sigma_6^{-1}$.
\end{itemize}
\end{lema}

\pf (i). First, \begin{align}\label{eq:sigma:tau} \sigma_1\tau_6=
\sigma_1 \, \sigma_2 \tau_3  \sigma_2^{-1}= \tau_3 \sigma_2
\tau_3^{-1 } \, \tau_3 \sigma_6 \tau_3^{-1 }  \,  \tau_3
\sigma_2^{-1}= \tau_3 \, \sigma_2 \sigma_6 \sigma_2^{-1}= \tau_3
\sigma_5= \sigma_5 \tau_3.
\end{align}
Applying now $\sigma_2 \,  \trid \underline{\quad}$ to
\eqref{eq:sigma:tau} we get $\sigma_3\tau_5=\tau_6\sigma_1$.
Applying  $\sigma_2 \, \trid \underline{\quad}$ to this last
identity, we have $\sigma_6\tau_1=\tau_5\sigma_3$. The rest is
similar.

(ii). By (i) and Lemma \ref{le:^4iguales} (ii) for $(\tau_i)_{1\leq
i \leq 6}$, we have that
\begin{align*}
\sigma_2^{-1}\tau_2=\sigma_2^{-1} \tau_4^{-1} \tau_4 \tau_2=
\sigma_1^{-1} \tau_6^{-1} \tau_1 \tau_6= \sigma_1^{-1} \tau_1.
\end{align*}
The other relations can be obtained in an analogous way.

(iii). It is easy to see that
\begin{align*}
\tau_2^{-2}\sigma_5\tau_5&= \tau_2^{-4}\tau_2\tau_2\tau_5\sigma_5=
\tau_1^{-4}\tau_2\tau_5\tau_1\sigma_5=
\tau_1^{-4}\tau_5\tau_1\tau_1\sigma_5\\&=
\tau_1^{-4}\tau_5\tau_1\sigma_4\tau_1=
\tau_1^{-4}\tau_5\sigma_3\tau_1\tau_1=
\tau_1^{-4}\tau_1\sigma_6\tau_1^2= \tau_1^{-1}\sigma_6.
\end{align*}

(iv) follows from (iii) applying $\sigma_2 \trid(\sigma_2 \trid
\underline{\quad})$.

(v). Clearly,
\begin{align*}
\sigma_2^{-2}\sigma_5\tau_5&=
\sigma_2^{-4}\sigma_2\sigma_2\sigma_5\tau_5=
\sigma_1^{-4}\sigma_2\sigma_5\sigma_1\tau_5=
\sigma_1^{-4}\sigma_5\sigma_1\sigma_1\tau_5=
\sigma_1^{-4}\sigma_5\sigma_1\tau_4\sigma_1\\&=
\sigma_1^{-4}\sigma_5\tau_3\sigma_1\sigma_1=
\sigma_1^{-4}\sigma_1\tau_6\sigma_1\sigma_1=\sigma_1^{-1}\tau_6=
\sigma_1^{-2}\tau_1\sigma_6.
\end{align*}

(vi) follows from (v) applying $\sigma_2 \trid(\sigma_2 \trid
\underline{\quad})$.\epf


\subsection{Applications}

Let $G$ be a finite group, $\oc$ a conjugacy class of $G$. Let
$(\sigma_i)_{1\leq i \leq 6}\subset \oc$ be of type $\oct$. We
define
\begin{align}
g_1:=\sigma_1,  \quad  g_2:=\sigma_5,  \quad   g_3:=\sigma_2,
\quad g_4:=\sigma_3,\quad g_5:=\sigma_4, \quad
g_6:=\sigma_2^2\sigma_1;
\end{align}
then, $\sigma_i=g_i \trid \sigma_1$, $1\leq i \leq 6$. It is easy
to see that following relations hold
\begin{align*}
\sigma_1 g_1&=g_1 \sigma_1,  &\sigma_1 g_2&=g_5 \sigma_1,
&\sigma_1 g_3&=g_2 \sigma_1,\\
\sigma_2 g_1&=g_3 \sigma_1,  &\sigma_2 g_2&=g_2\sigma_1, &\sigma_2
g_3&=g_6 \sigma_1^{-1},\\
\sigma_3 g_1&=g_4 \sigma_1,  &\sigma_3 g_2&=g_1\sigma_6,
&\sigma_3g_3&=g_3 \sigma_1,\\
\sigma_4 g_1&=g_5 \sigma_1,  &\sigma_4 g_2&=g_2\sigma_6,
&\sigma_4g_3&=g_1 \sigma_6,\\
\sigma_5g_1&=g_2 \sigma_1,  &\sigma_5
g_2&=g_6\sigma_1^{-2}\sigma_6,  &\sigma_5g_3&=g_3 \sigma_6,\\
\sigma_6g_1&=g_1 \sigma_6,  &\sigma_6 g_2&=g_3\sigma_6,
&\sigma_6g_3&=g_4 \sigma_6,
\end{align*}
\begin{align*}
\sigma_1 g_4&=g_3 \sigma_1, &\sigma_1 g_5&=g_4 \sigma_1, &\sigma_1
g_6&=g_6 \sigma_6,\\
\sigma_2 g_4&=g_4 \sigma_6,  &\sigma_2 g_5&=g_1 \sigma_6,
&\sigma_2 g_6&=g_5 \sigma_6^3,\\
\sigma_3 g_4&=g_6 \sigma_6^{-1}, &\sigma_3 g_5&=g_5 \sigma_6,
&\sigma_3g_6&=g_2 \sigma_1^3,\\
\sigma_4 g_4&=g_4 \sigma_1, &\sigma_4 g_5&=g_6
\sigma_1\sigma_6^{-2},&\sigma_4 g_6&=g_3 \sigma_1^{2}\sigma_6,\\
\sigma_5 g_4&=g_1 \sigma_6,  &\sigma_5
g_5&=g_5\sigma_1,  &\sigma_5 g_6&=g_4 \sigma_1 \sigma_6^2,\\
\sigma_6 g_4&=g_5 \sigma_6, &\sigma_6 g_5&=g_2 \sigma_6, &\sigma_6
g_6&=g_6 \sigma_1.
\end{align*}

Let $\rho=(\rho,V) \in \widehat{G^{\sigma_1}}$ and $v\in V-0$. Assume that
$v$ is an eigenvector of $\rho(\sigma_6)$ with eigenvalue $\lambda$.
We
define $W:=$ span- $\{g_iv \, | \, 1\leq i \leq 6\}$. Then, $W$ is
a braided vector subspace of $M(\oc, \rho)$.

\begin{lema}\label{le:4-ciclos}
Let $(\sigma_i)_{1\leq i \leq 6}$, $(g_i)_{1\leq i \leq 6}$,
$(\rho,V)\in \widehat{G^{\sigma_1}}$, $W$, $\lambda$ as above. Assume that
$q_{\sigma_1\sigma_1}= \lambda =-1$. Then $W\simeq M(\oc_4^4, \chi_-)$ as
braided vector spaces.
\end{lema}
\pf Since $q_{\sigma_1\sigma_1}=-1$ we have that
$\rho(\sigma_i^4)=\Id$, $1\leq i \leq 6$, from Lemma
\eqref{le:^4iguales} (i). Let $\widetilde{\sigma}_i$ be as in
\eqref{eq:relss4}. If we choose
\begin{align*}
\widetilde{g}_1=\widetilde{\sigma}_1, \quad
\widetilde{g}_2=\widetilde{\sigma}_5, \quad
\widetilde{g}_3=\widetilde{\sigma}_2, \quad
\widetilde{g}_4=\widetilde{\sigma}_3, \quad
\widetilde{g}_5=\widetilde{\sigma}_4, \quad
\widetilde{g}_6=\widetilde{\sigma}_2^2\widetilde{\sigma}_1,
\end{align*}
then $\widetilde{g}_i\trid
\widetilde{\sigma}_1=\widetilde{\sigma}_i$, $1\leq i \leq 6$.
Thus, $M(\oc_4^4, \chi_-)=$ span-$\{\widetilde{g}_iv_0,\, | \,
1\leq i \leq 6 \}$, with $v_0\in V_0-0$, where $V_0$ is the vector
space affording the representation $\chi_-$ of
$\s_4^{(1\,2\,3\,4)}$. Now, the map $W\to M(\oc_4^4, \chi_-)$
given by $g_iv \mapsto \widetilde{g}_iv_0$, $1\leq i \leq 6$, is
an isomorphism of braided vector spaces. \epf

The next lemma is needed for the main result of the section.

\begin{lema}\label{le:rel:octadoble}
Let $\sigma_i$, $\tau_i$, $1\leq i \leq 6$, be distinct elements
in $G$, $\oc$ a conjugacy class of $G$. Assume that
$(\sigma,\tau)\subseteq \oc$ is of type $\oct^{(2)}$, with $g\in
G$ such that $g\trid \sigma_1 =\tau_1$. Let
\begin{equation}\label{eq:losgparaR4^2}
\begin{aligned}
g_1&:=\sigma_1, \quad &g_2&:=\sigma_5, \quad &g_3&:=\sigma_2,
\quad &g_4&:=\sigma_3, \\ g_5&:=\sigma_4, \quad
&g_6&:=\sigma_2^2\sigma_1, &g_7&:=g\sigma_1, \quad &g_8&:=\tau_5g, \\
g_9&:=\tau_2g, \quad &g_{10}&:=\tau_3g, \quad &g_{11}&:=\tau_4g,
\quad &g_{12}&:=\tau_2^2g\sigma_1.
\end{aligned}
\end{equation}

Then, the following relations hold:
\begin{align*}
\tau_1 g_7&=g_7 \sigma_1,  &\tau_1 g_8&=g_{11} \sigma_1, &\tau_1
g_9&=g_8 \sigma_1,\\
\tau_2 g_7&=g_9 \sigma_1,  &\tau_2 g_8&=g_8\sigma_1, &\tau_2
g_9&=g_{12} \sigma_1^{-1},\\
\tau_3 g_7&=g_{10} \sigma_1,  &\tau_3 g_8&=g_7g^{-1}\tau_6 g,
&\tau_3g_9&=g_9 \sigma_1,\\
\tau_4 g_7&=g_{11} \sigma_1,  &\tau_4 g_8&=g_8 g^{-1}\tau_6 g,
&\tau_4g_9&=g_7 g^{-1}\tau_6 g,\\
\tau_5g_7&=g_8 \sigma_1,  &\tau_5
g_8&=g_{12}\sigma_1^{-2}g^{-1}\tau_6 g, &\tau_5g_9&=g_9
g^{-1}\tau_6 g,\\
\tau_6g_7&=g_7 g^{-1}\tau_6 g,  &\tau_6 g_8&=g_9g^{-1}\tau_6 g,
&\tau_6g_9&=g_{10} g^{-1}\tau_6 g,
\end{align*}
\begin{align*}
\tau_1 g_{10}&=g_9 \sigma_1, &\tau_1 g_{11}&=g_{10} \sigma_1,
&\tau_1g_{12}&=g_{12} g^{-1}\tau_6 g,\\
\tau_2 g_{10}&=g_{10} g^{-1}\tau_6 g, &\tau_2 g_{11}&=g_7
g^{-1}\tau_6 g, &\tau_2 g_{12}&=g_{11} (g^{-1}\tau_6 g)^3,\\
\tau_3 g_{10}&=g_{12} (g^{-1}\tau_6 g)^{-1}, &\tau_3
g_{11}&=g_{11} g^{-1}\tau_6 g, &\tau_3
g_{12}&=g_8 \sigma_1^3,\\
\tau_4 g_{10}&=g_{10} \sigma_1, &\tau_4 g_{11}&=g_{12}
\sigma_1(g^{-1}\tau_6 g)^{-2},
&\tau_4 g_{12}&=g_9 \sigma_1^{2}g^{-1}\tau_6 g,\\
\tau_5 g_{10}&=g_7 g^{-1}\tau_6 g, &\tau_5
g_{11}&=g_{11} \sigma_1,  &\tau_5 g_{12}&=g_{10} \sigma_1 (g^{-1}\tau_6 g)^2,\\
\tau_6 g_{10}&=g_{11} g^{-1}\tau_6 g, &\tau_6 g_{11}&=g_8
g^{-1}\tau_6 g, &\tau_6 g_{12}&=g_{12} \sigma_1,
\end{align*}
\begin{align*}
\sigma_1 g_7&=g_7 g^{-1}\sigma_1 g,  &\sigma_1 g_8&=g_{11}
g^{-1}\sigma_1 g, &\sigma_1 g_9&=g_8 g^{-1}\sigma_1 g,\\
\sigma_2 g_7&=g_9 g^{-1}\sigma_1 g,  &\sigma_2
g_8&=g_8g^{-1}\sigma_1 g, &\sigma_2 g_9&=g_{12}
\sigma_1^{-2}(g^{-1}\sigma_1 g),\\
\sigma_3 g_7&=g_{10} g^{-1}\sigma_1 g,  &\sigma_3
g_8&=g_7g^{-1}\sigma_6 g, &\sigma_3g_9&=g_9 g^{-1}\sigma_1 g,\\
\sigma_4 g_7&=g_{11} g^{-1}\sigma_1 g,  &\sigma_4 g_8&=g_8
g^{-1}\sigma_6 g, &\sigma_4g_9&=g_7 g^{-1}\sigma_6 g,\\
\sigma_5g_7&=g_8 g^{-1}\sigma_1 g,  &\sigma_5
g_8&=g_{12}\sigma_1^{-2}g^{-1}\sigma_6 g, &\sigma_5g_9&=g_9
g^{-1}\sigma_6 g,\\
\sigma_6g_7&=g_7 g^{-1}\sigma_6 g,  &\sigma_6
g_8&=g_9g^{-1}\sigma_6 g, &\sigma_6g_9&=g_{10} g^{-1}\sigma_6 g,
\end{align*}
\begin{align*}
\sigma_1 g_{10}&=g_9 g^{-1}\sigma_1 g, &\sigma_1 g_{11}&=g_{10}
g^{-1}\sigma_1 g, &\sigma_1
g_{12}&=g_{12} g^{-1}\sigma_6 g,\\
\sigma_2 g_{10}&=g_{10} g^{-1}\sigma_6 g, &\sigma_2 g_{11}&=g_7
g^{-1}\sigma_6 g, &\sigma_2g_{12}&=g_{11} \gamma_{2,12},\\
\sigma_3 g_{10}&=g_{12} \gamma_{3,10}, &\sigma_3 g_{11}&=g_{11}
g^{-1}\sigma_6 g, &\sigma_3g_{12}&=g_8 \sigma_1^2(g^{-1}\sigma_1 g),\\
\sigma_4 g_{10}&=g_{10} g^{-1}\sigma_1 g, &\sigma_4 g_{11}&=g_{12}
\gamma_{4,11},
&\sigma_4 g_{12}&=g_9 \sigma_1^{2}g^{-1}\sigma_6 g,\\
\sigma_5 g_{10}&=g_7 g^{-1}\sigma_6 g, &\sigma_5 g_{11}&=g_{11}
g^{-1}\sigma_1 g,
&\sigma_5 g_{12}&=g_{10}\gamma_{5,12} ,\\
\sigma_6 g_{10}&=g_{11} g^{-1}\sigma_6 g, &\sigma_6 g_{11}&=g_8
g^{-1}\sigma_6 g, &\sigma_6 g_{12}&=g_{12} g^{-1}\sigma_1 g,
\end{align*}
where $\gamma_{2,12}= \sigma_1^2(g^{-1}\sigma_1
g)^{-2}(g^{-1}\sigma_6 g)^3$, $\gamma_{3,10}=
\sigma_1^{-2}(g^{-1}\sigma_1 g)^{2}(g^{-1}\sigma_6 g)^{-1}$,
$\gamma_{4,11}=\sigma_1^{-2}(g^{-1}\sigma_1 g)^{3}(g^{-1}\sigma_6
g)^{-2}$ and $\gamma_{5,12}=\sigma_1^2 (g^{-1}\sigma_1
g)^{-1}(g^{-1}\sigma_6 g)^2$,

\begin{align*}
\tau_1 g_1&=g_1 \tau_1,  &\tau_1 g_2&=g_5 \tau_1, &\tau_1 g_3&=g_2
\tau_1,\\
\tau_2 g_1&=g_3 \tau_1,  &\tau_2 g_2&=g_2\tau_1, &\tau_2 g_3&=g_6
\sigma_1^{-2}\tau_1,\\
\tau_3 g_1&=g_4 \tau_1,  &\tau_3 g_2&=g_1\tau_6, &\tau_3 g_3&=g_3
\tau_1,\\
\tau_4 g_1&=g_5 \tau_1,  &\tau_4 g_2&=g_2\tau_6, &\tau_4 g_3&=g_1
\tau_6,\\
\tau_5 g_1&=g_2 \tau_1,  &\tau_5 g_2&=g_6\sigma_1^{-2}\tau_6,
&\tau_5 g_3&=g_3 \tau_6,\\
\tau_6 g_1&=g_1 \tau_6,  &\tau_6 g_2&=g_3\tau_6, &\tau_6 g_3&=g_4
\tau_6,
\end{align*}
\begin{align*}
\tau_1 g_4&=g_3 \tau_1, &\tau_1 g_5&=g_4
\tau_1,  &\tau_1 g_6&=g_6 \tau_6,\\
\tau_2 g_4&=g_4 \tau_6,  &\tau_2
g_5&=g_1 \tau_6,  &\tau_2 g_6&=g_5 \sigma_1^3\tau_1\sigma_6,\\
\tau_3 g_4&=g_6 \sigma_1^{-1}\tau_1\sigma_6^{-1},
&\tau_3 g_5&=g_5 \tau_6,  &\tau_3 g_6&=g_2 \sigma_1^2\tau_1,\\
\tau_4 g_4&=g_4 \tau_1, &\tau_4
g_5&=g_6 \tau_1\sigma_6^{-2}, &\tau_4 g_6&=g_3 \sigma_1\tau_1\sigma_6,\\
\tau_5 g_4&=g_1 \tau_6,  &\tau_5
g_5&=g_5 \tau_1,  &\tau_5 g_6&=g_4 \tau_1 \sigma_6^2,\\
\tau_6 g_4&=g_5 \tau_6, &\tau_6 g_5&=g_2 \tau_6, &\tau_6 g_6&=g_6
\tau_1.
\end{align*}
\end{lema}

\pf The proof follows by straightforward computations, Lemma
\ref{le:^4iguales} for $\sigma$ and $\tau$, and Lemma
\ref{le:sigma:tau}.\epf

Here is the main result of this section.

\begin{theorem}\label{teor:aplicAHSch2}
Let $\sigma_i$, $\tau_i\in G$, $1\leq i \leq 6$, distinct elements
in $G$, $\oc$ a conjugacy class of $G$ and $\rho=(\rho,V)\in
\widehat{G^{\sigma_1}}$. Let us suppose that
\begin{itemize}
\item[(H1)] $(\sigma,\tau)$ is of type $\oct^{(2)}$,
\item[(H2)] $(\sigma,\tau)\subseteq \oc$, with
$g\in G$ such that $g\trid \sigma_1 =\tau_1$,
\item[(H3)] $q_{\sigma_1\sigma_1}=-1$,
\end{itemize}
there exists $v\in V-0$ such that
\begin{itemize}
\item[(H4)] $\rho(\sigma_6) v=-v$,
\item[(H5)] $\rho(\tau_1) v=-v$,
\end{itemize}
and there exists $w\in V-0$ such that
\begin{itemize}
\item[(H6)] $\rho(g^{-1}\sigma_1g) w=-w$,
\item[(H7)] $\rho(g^{-1}\sigma_6g) w=-w$,
\end{itemize}
Then $\dim\toba(\oc,\rho)=\infty$.
\end{theorem}

\pf Let $g_j\in G$ , $1\leq j \leq 12$, as in
\eqref{eq:losgparaR4^2}. Then, $g_j\trid \sigma_1=\sigma_j$,
$1\leq j \leq 6$, and $g_j\trid \sigma_1=\tau_{j-6}$, $7\leq j \leq
12$. By Lemma \ref{le:rel:octadoble}, we have that
\begin{itemize}
\item[(a)] if $1\leq i,j \leq 6$, then $g_{i\trid j}^{-1} \sigma_i
g_j=\sigma_1^r\sigma_6^s$, with $r+s$ odd,
\item[(b)] if $7\leq i,j \leq 12$, then $g_{i\trid j}^{-1} \tau_{i-6} g_j=
\sigma_1^r(g^{-1} \tau_6 g)^s$, with $r+s$ odd,
\item[(c)] if $1\leq i \leq 6$ and $7\leq j \leq 12$, then
$g_{i\trid j}^{-1} \sigma_{i} g_j= \sigma_1^r(g^{-1}\sigma_1 g)^s
(g^{-1}\sigma_6 g)^t$, with $r+s+t$ odd,
\item[(d)] if $1\leq j \leq 6$ and $7\leq i \leq 12$, then $g_{i\trid j}^{-1} \tau_{i-6} g_j=
\sigma_1^r\tau_1^s\sigma_6^t$, with $r+s+t$ odd, because
$\tau_6=\sigma_1^{-1}\tau_1\sigma_6$.
\end{itemize}

Let $W:=$ span-$\{g_iv, \, | \, 1\leq i \leq 6\}$ and $W':=$
span-$\{g_iw,\, | \, 7\leq i \leq 12\}$, with $v$, $w\in V-0$,
where $v$ satisfies (H4)-(H5) and $w$ satisfies (H6)-(H7). Then,
$W$ and $W'$ are braided vector subspaces of $M(\oc,\rho)$. We
will prove that
\begin{align*}
W\oplus W' \simeq M(\oc_4^4, \chi_-) \oplus M(\oc_4^4, \chi_-),
\end{align*}
as braided vector spaces. Hence $\dim\toba(W\oplus W')=\infty$, by
Theorem \ref{teor:AHS2}, and the result follows from Lemma
\ref{trivialbraiding}.

By Remark \ref{obs:isoEVTocta}, we only need to see that the
isomorphism of linear vector spaces $W\oplus W' \to M(\oc_4^4, \chi_-) \oplus M(\oc_4^4, \chi_-)$ given by
\begin{align*}
g_iv \mapsto \widetilde g_i \quad \text{and} \quad g_{i+6}w \mapsto \widetilde g_{i+6} \qquad 1\leq i \leq 6,
\end{align*}
respects the braiding, and this is just a matter of the cocycle. For this, we
compute explicitly the braiding in the basis $\{g_iv, \, g_{j+6}w,
\, | \, 1\leq i,j \leq 6\}$ of $W\oplus W'$.


By (a), (H3) and (H4), if $1\leq i,j \leq 6$, then
\begin{align*}
c(g_i v\ot g_j v)=g_{i\trid j} \rho(g_{i\trid j}^{-1} \sigma_i
g_j)(v) \ot g_i v=- g_{i\trid j} v\ot g_iv.
\end{align*}

\noindent From Lemma \ref{le:sigma:tau} (i),
$\tau_6=\sigma_1^{-1}\tau_1\sigma_6$. Thus,
$g^{-1}\tau_6g=(g^{-1}\sigma_1g)^{-1}\sigma_1(g^{-1}\sigma_6g)$.
By (b), (H3), (H6) and (H7), if $7\leq i,j \leq 12$, then
\begin{align*}
c(g_i w\ot g_j w)=g_{i\trid j} \rho(g_{i\trid j}^{-1} \tau_{i-6}
g_j)(w) \ot g_i w=- g_{i\trid j} w\ot g_iw.
\end{align*}

\noindent By (c), (H3), (H6) and (H7), if $1\leq i \leq 6$ and
$7\leq j \leq 12$, then
\begin{align*}
c(g_i v\ot g_j w)=g_{i\trid j} \rho(g_{i\trid j}^{-1} \sigma_i
g_j)(w) \ot g_i v=- g_{i\trid j} w\ot g_iv.
\end{align*}

\noindent By (d), (H3), (H4) and (H5), if $1\leq j \leq 6$ and
$7\leq i \leq 12$, then
\begin{align*}
c(g_i w\ot g_j v)=g_{i\trid j} \rho(g_{i\trid j}^{-1} \tau_{i-6}
g_j)(v) \ot g_i w=- g_{i\trid j} v\ot g_iw.
\end{align*}
This completes the proof.\epf

As an immediate consequence we have the following result.

\begin{cor}\label{co:especial2}
Let $\sigma_i$, $\tau_i\in G$, $1\leq i \leq 6$ all distinct,
$\oc$ a conjugacy class of $G$ and $\rho=(\rho,V)\in
\widehat{G^{\sigma_1}}$ with $q_{\sigma_1\sigma_1}=-1$. Assume
that $(\sigma,\tau)\subseteq \oc$ is of type $\oct^{(2)}$. If
$\sigma_6=\sigma_1^d$ and $\tau_1=\sigma_1^e$ for some $d$, $e \in
\Z$, then $\dim\toba(\oc,\rho)=\infty$.
\end{cor}

\pf Note that $d$ and $e$ are odd, since they are relatively prime
with $\vert\sigma_1\vert$. Hence the hypothesis (H4) and (H5)
hold. Now $g^{-1}\sigma_1 g=\sigma_1^{e^{|g|-1}}$. Then
$\rho(g^{-1}\sigma_1 g)=-\Id$ and (H6) holds. The proof of (H7) is
similar. \epf

\begin{exa}\label{ex:8-ciclo}
Let $m\geq 8$. Let $\sigma\in \sm$ of type
$(1^{n_1},2^{n_2},8^{n_8})$, with $n_8\geq 1$, $\oc$ the conjugacy
class of $\sigma$ and $\rho \in \widehat{\s_m^{\sigma}}$.  Then
$\dim \toba(\oc,\rho)=\infty$.
\end{exa}

\pf By Lemma \ref{odd}, we may suppose that $q_{\sigma\sigma}=
-1$. If $n_8\geq 3$, then $\dim \toba(\oc,\rho)=\infty$, from
Corollary \ref{co:type2k^3}. We consider two cases.

\emph{CASE (I):} $n_8=1$. Let
$A_{8}=(i_1\,\,i_2\,\,i_3\,\,i_4\,\,i_5\,\,i_6\,\,i_7\,\,i_8)$ the
8-cycle appearing in the decomposition of $\sigma$ as product of
disjoint cycles. We set $\alpha:=\sigma  \, A_8^{-1}$ and define
$\sigma_1:=\sigma$, $\sigma_6:=\sigma_1^3$, $\tau_1:=\sigma_1^5$,
$\tau_6:=\sigma_1^{-1}$,
\begin{align*}
\sigma_2&:=(i_1\,\,i_3\,\,i_8\,\,i_6\,\,i_5\,\,i_7\,\,i_4\,\,i_2)\,\alpha,
&\sigma_3&:=(i_1\,\,i_8\,\,i_2\,\,i_7\,\,i_5\,\,i_4\,\,i_6\,\,i_3)\,\alpha,\\
\sigma_4&:=(i_1\,\,i_6\,\,i_4\,\,i_3\,\,i_5\,\,i_2\,\,i_8\,\,i_7)\,\alpha,
&\sigma_5&:=(i_1\,\,i_7\,\,i_6\,\,i_8\,\,i_5\,\,i_3\,\,i_2\,\,i_4)\,\alpha,\\
\tau_2&:=(i_1\,\,i_7\,\,i_8\,\,i_2\,\,i_5\,\,i_3\,\,i_4\,\,i_6)\,\alpha,
&\tau_3&:=(i_1\,\,i_4\,\,i_2\,\,i_3\,\,i_5\,\,i_8\,\,i_6\,\,i_7)\,\alpha,\\
\tau_4&:=(i_1\,\,i_2\,\,i_4\,\,i_7\,\,i_5\,\,i_6\,\,i_8\,\,i_3)\,\alpha,
&\tau_5&:=(i_1\,\,i_3\,\,i_6\,\,i_4\,\,i_5\,\,i_7\,\,i_2\,\,i_8)\,\alpha.
\end{align*}

\emph{CASE (II):} $n_8=2$. Let
\begin{align*}
A_{1,8}=(i_1\,\,i_2\,\,i_3\,\,i_4\,\,i_5\,\,i_6\,\,i_7\,\,i_8)
\,\, \text{ and } \,\,
A_{2,8}=(i_9\,\,i_{10}\,\,i_{11}\,\,i_{12}\,\,i_{13}\,\,i_{14}\,\,i_{15}\,\,i_{16})
\end{align*}
the two 8-cycles appearing in the decomposition of $\sigma$ as
product of disjoint cycles. We call $A_8=A_{1,8}A_{2,8}$,
$\alpha:=\sigma  \, A_8^{-1}$ and define $\sigma_1:=\sigma$,
$\sigma_6:=\sigma_1^3$, $\tau_1:=\sigma_1^5$,
$\tau_6:=\sigma_1^{-1}$,
\begin{align*}
\sigma_2&:=(i_1\,\,i_3\,\,i_8\,\,i_6\,\,i_5\,\,i_7\,\,i_4\,\,i_2)(i_{9}\,\,i_{11}\,\,i_{16}\,\,i_{14}\,\,i_{13}\,\,i_{15}\,\,i_{12}\,\,i_{10})\,\alpha,\\
\sigma_3&:=(i_1\,\,i_8\,\,i_2\,\,i_7\,\,i_5\,\,i_4\,\,i_6\,\,i_3)(i_{9}\,\,i_{16}\,\,i_{10}\,\,i_{15}\,\,i_{13}\,\,i_{12}\,\,i_{14}\,\,i_{11})\,\alpha,\\
\sigma_4&:=(i_1\,\,i_6\,\,i_4\,\,i_3\,\,i_5\,\,i_2\,\,i_8\,\,i_7)(i_{9}\,\,i_{14}\,\,i_{12}\,\,i_{11}\,\,i_{13}\,\,i_{10}\,\,i_{16}\,\,i_{15})\,\alpha,\\
\sigma_5&:=(i_1\,\,i_7\,\,i_6\,\,i_8\,\,i_5\,\,i_3\,\,i_2\,\,i_4)(i_{9}\,\,i_{15}\,\,i_{14}\,\,i_{16}\,\,i_{13}\,\,i_{11}\,\,i_{10}\,\,i_{12})\,\alpha,\\
\tau_2&:=(i_1\,\,i_7\,\,i_8\,\,i_2\,\,i_5\,\,i_3\,\,i_4\,\,i_6)(i_{9}\,\,i_{15}\,\,i_{16}\,\,i_{10}\,\,i_{13}\,\,i_{11}\,\,i_{12}\,\,i_{14})\,\alpha,\\
\tau_3&:=(i_1\,\,i_4\,\,i_2\,\,i_3\,\,i_5\,\,i_8\,\,i_6\,\,i_7)(i_{9}\,\,i_{12}\,\,i_{10}\,\,i_{11}\,\,i_{13}\,\,i_{16}\,\,i_{14}\,\,i_{15})\,\alpha,\\
\tau_4&:=(i_1\,\,i_2\,\,i_4\,\,i_7\,\,i_5\,\,i_6\,\,i_8\,\,i_3)(i_{9}\,\,i_{10}\,\,i_{12}\,\,i_{15}\,\,i_{13}\,\,i_{14}\,\,i_{16}\,\,i_{11})\,\alpha,\\
\tau_5&:=(i_1\,\,i_3\,\,i_6\,\,i_4\,\,i_5\,\,i_7\,\,i_2\,\,i_8)(i_{9}\,\,i_{11}\,\,i_{14}\,\,i_{12}\,\,i_{13}\,\,i_{15}\,\,i_{10}\,\,i_{16})\,\alpha.
\end{align*}

In both cases, $\sigma_6=\sigma_1^3$ and $\tau_1=\sigma_1^5$ and
$(\sigma,\tau)\subseteq \oc$ is of type $\oct^{(2)}$. Then the
result follows from Corollary \ref{co:especial2}.\epf

\begin{rmk}
(i). The discussion in the preceding example can be adapted to
$\sigma\in \sm$ of type $(1^{n_1},2^{n_2},\dots,m^{n_m})$ provided
that $n_8 \geq 1$; but then some requirements on the
representation $\rho$ have to be imposed.

(ii). Let $N = 2^n$ with $n\ge 4$. It can be shown that the orbit
of the $N$-cycle in $\s_{N}$ contains no family of type $\oct$
using Lemma \ref{le:^4iguales}.

(iii). The orbit with label $j=4$ of the Mathieu group $M_{22}$
contains a family of type $\oct^{(2)}$, and therefore this group
admits no finite-dimensional pointed Hopf algebra except the group
algebra itself \cite{F}.
\end{rmk}

\subsection*{Acknowledgement}
The authors are grateful to the referee for carefully reading the
paper and for his/her comments.


\begin{thebibliography}{AFGV}

\bibitem[AF1]{AF}
N. Andruskiewitsch and F. Fantino, \emph{On pointed Hopf algebras
associated with unmixed conjugacy classes in $\mathbb S_m$}, J.
Math. Phys. \textbf{48} (2007), 033502, 1-26.

\bibitem[AF2]{AF2} \bysame,
\emph{On pointed Hopf algebras associated with alternating and
dihedral groups}, Rev. Uni\'on Mat. Argent. 48-3, (2007), 57-71,

\bibitem[AFGV]{AFV}
N. Andruskiewitsch,  F. Fantino, M. Gra\~na and  L. Vendramin,
\emph{On pointed Hopf algebras associated to sporadic groups}, in
preparation.

\bibitem[AFZ]{AFZ}
N. Andruskiewitsch,  F. Fantino and   S. Zhang, \emph{On pointed
Hopf algebras associated to symmetric groups},
\texttt{arXiv:0807.2406v2 [math.QA]}.

\bibitem[AG]{AG1} N. Andruskiewitsch and  M. Gra\~na,
\emph{From racks to pointed Hopf algebras}, Adv. Math.
\textbf{178}  (2003), 177 -- 243.

\bibitem[AHS]{AHS} N. Andruskiewitsch, I. Heckenberger and H.-J. Schneider,
{\em The Nichols algebra of a semisimple Yetter-Drinfeld module},
\texttt{arXiv:0803.2430v1 [math.QA]}.

\bibitem[AS1]{AS1}  N. Andruskiewitsch and  H.-J. Schneider, \textit{Lifting of
quantum linear spaces and pointed Hopf algebras of order $p^3$}.
J. Algebra 209, 658-691 (1998).

\bibitem[AS2]{AS-adv}  \bysame, {\em Finite quantum groups and Cartan matrices}, Adv.
Math. \textbf{154} (2000), 1--45.

\bibitem[AS3]{AS-cambr} \bysame, \emph{Pointed Hopf Algebras}, in ``New directions in Hopf
algebras'', 1--68, Math. Sci. Res. Inst. Publ. \textbf{43},
Cambridge Univ. Press, Cambridge, 2002.

\bibitem[AS4]{AS6}\bysame, \textit{On the
classification of finite-dimensional pointed Hopf algebras}. Ann.
Math., accepted, 43 pp., preprint \texttt{math.QA/0502157}.

\bibitem[AZ]{AZ}  N. Andruskiewitsch and  S. Zhang, \emph{On pointed Hopf algebras associated to some conjugacy
classes in $\mathbb S_n$}, Proc. Amer. Math. Soc. 135  (2007),
2723-2731.


\bibitem[F1]{F}
 F. Fantino, \emph{On pointed Hopf algebras associated with
Mathieu groups},  \texttt{arXiv: 0711.3142v2 [math.QA]}. 

\bibitem[F2]{F-tesis} \bysame, \emph{\'Algebras de Hopf punteadas
sobre grupos no abelianos}. Tesis de doc\-to\-ra\-do, Universidad
Nacional de C\'ordoba (2008).
\texttt{www.mate.uncor.edu/\~{}fantino/}.


\bibitem[FGV]{FGV}  S. Freyre,  M. Gra\~na and  L. Vendramin,
\emph{On Nichols algebras over $\mathbf{GL}(2,\mathbb{F}_q)$ and
${\mathbf{SL}(2,\mathbb{F}_q)}$}, J. Math. Phys. \textbf{48}
(2007), 123513-1 -- 123513-11.

\bibitem[FV]{FV}  S. Freyre and  L. Vendramin,
\emph{On Nichols algebras over $\mathbf{PSL}(2,\mathbb{F}_q)$ and \newline
${\mathbf{PGL}(2,\mathbb{F}_q)}$}, \texttt{arXiv:0802.2567v1
[math.QA]}.

\bibitem[FK]{FK} S. Fomin and  K. N. Kirillov,  {\it Quadratic algebras, Dunkl
elements, and Schubert calculus}, Progr. Math.    {\bf 172},
Birkhauser, (1999),  146--182.




\bibitem[G1]{G1}  M. Gra\~na, \emph{On Nichols algebras of low dimension},
 Contemp. Math. {\bf 267}  (2000),  111--134.

\bibitem[G2]{G2}\bysame, Finite dimensional Nichols algebras of non-diagonal group
type, zoo of examples available at
\texttt{http://mate.dm.uba.ar/~matiasg/zoo.html}.

\bibitem[H1]{H-inv}  I. Heckenberger,
{\em The Weyl groupoid of a Nichols algebra of diagonal type},
Invent. Math. 164 (2006), 175--188.

\bibitem[H2]{H-all}
\bysame, \textit{Classification of arithmetic root systems}, Adv.
in Math., accepted, preprint \texttt{math.QA/0605795}.

\bibitem[MS]{MS}
 A. Milinski and H-J.  H.-J. Schneider, \emph{Pointed
Indecomposable Hopf Algebras over Coxeter Groups}, Contemp. Math.
\textbf{267} (2000), 215--236.


\bibitem[S]{S} Jean-Pierre Serre, \emph{Linear representations of finite groups},
Springer-Verlag, New York 1977.



\end{thebibliography}
\end{document}